\documentclass{amsart}
\usepackage{amssymb}
\usepackage{scalerel}
\usepackage{comment}

\newtheorem{thm}{Theorem}[section]
\newtheorem{prop}[thm]{Proposition}
\newtheorem{lem}[thm]{Lemma}
\newtheorem{cor}[thm]{Corollary}

\theoremstyle{definition}
\newtheorem{definition}[thm]{Definition}

\theoremstyle{remark}

\newtheorem{remark}[thm]{Remark}

\numberwithin{equation}{section}

\newcommand{\ip}[2]{\left\langle {#1} , {#2} \right\rangle} 
\newcommand{\norm}[1]{\left\lVert #1 \right\rVert} 
\newcommand{\R}{\mathbb{R}} 
\newcommand{\C}{\mathbb{C}} 
\newcommand{\T}{\mathbb{T}} 

\newcommand{\tr}{\text{tr}} 
\newcommand{\Falpha}{\mathcal{F}^{2}_{\alpha}}
\newcommand{\calF}{\mathcal{F}}
\newcommand{\Weyl}{W}

\begin{document}


\title{A Semi-Classical Szeg\H{o}-type Limit Theorem for Toeplitz Operators}


\author{Trevor Camper}
\address{Department of Mathematics, Dartmouth College,
Hanover, NH 03755}
\email{trevor.camper@dartmouth.edu}


\author{Mishko Mitkovski}\thanks{Supported in part by NSF Grant DMS-2453810.}
\address{School of Mathematical and Statistical Sciences, Clemson University,
Clemson, SC 29630}
\email{mmitkov@clemson.edu}



\begin{abstract}
We establish an abstract Szeg\H{o}-type limit theorem for Toeplitz operators associated with normalized continuous Parseval frames. Our assumptions consist essentially of a natural normalization of the frame measures and an approximate-identity condition for the normalized frame vectors. Under these hypotheses, the Szeg\H{o} asymptotics follow from a general argument based on Berezin--Lieb inequalities and
convergence of Berezin transforms, independently of the particular geometry of the underlying space.

As a principal application, we obtain a Szeg\H{o}-type limit theorem for Toeplitz operators on weighted Bergman spaces of the unit ball. This setting lies outside the standard Widom-type localization approach because the underlying hyperbolic measure has exponential volume growth. We also obtain Szeg\H{o} limit theorems for metrizable compact and locally compact Abelian groups, extending the classical compact-group results to more general symbols and recovering recent results for multidimensional tori as special cases. Further
applications include Gabor localization operators, Paley--Wiener spaces, eigenvalue distribution, and sampling density results.

\end{abstract}


\maketitle






\section{Introduction}

The classical Szeg\H o limit theorem \cite{SZ15}, sometimes also called the first Szeg\H o theorem, describes the asymptotic behavior of the spectra of Toeplitz matrices generated by a fixed symbol as the matrix size tends to infinity. More precisely, let $\sigma:\mathbb{T} \to \mathbb{R}$ be a symbol, and consider the compressions $P_n M_\sigma P_n : L^2(\mathbb{T}) \to L^2(\mathbb{T})$ of the multiplication operator $M_\sigma : L^2(\mathbb{T}) \to L^2(\mathbb{T})$, $M_\sigma f = \sigma f$, onto the span of the complex exponentials $\operatorname{span}\{ e^{ikx} : k = 0, \pm 1, \pm 2, \dots, \pm n \} \subseteq L^2(\mathbb{T})$. The matrices of these compressions, represented in this exponential basis, are the Toeplitz matrices $[\hat{\sigma}(i-j)]_{i,j=-n,\dots,n}$.

 One of the most classical forms of the Szeg\H o limit theorem states that for a strictly positive symbol $\sigma\in L^{1}(\mathbb{T})$ we have 

\begin{equation}\label{Szegolog}
\lim_{n\to\infty}\frac{1}{2n+1}\text{tr}\left(\log\left(P_{n}M_{\sigma}P_{n}\right)\right)=\frac{1}{2\pi}\int_{0}^{2\pi}\log\sigma\left(e^{i\theta}\right)d\theta. \end{equation}

\noindent Another, equally classical, version of the Szeg\H{o} limit theorem requires the symbol to be continuous, but not necessarily positive. In this case, for any function $\psi$ continuous on a closed interval containing the range of $\sigma$, we have 
\begin{equation}\label{Szegocont}
    \lim_{n\to\infty}\frac{1}{2n+1}\text{tr}\left(\psi\left(P_{n}M_{\sigma}P_{n}\right)\right)=\frac{1}{2\pi}\int_{0}^{2\pi}\psi\left(\sigma\left(e^{i\theta}\right)\right)d\theta. 
\end{equation}
Here, as everywhere above, $P_n$ denotes the orthogonal projection onto 
$\operatorname{span} \{ e^{ikx} : k = 0, \pm 1, \pm 2, \dots, \pm n \}$.

Numerous authors have explored extensions and generalizations of the first Szeg\H{o} limit theorem in various contexts and settings. A few notable works in this area are~\cite{guille79,kryein1983some,laptev1996szego,widom1979eigenvalue,widom2006asymptotic}. For a more comprehensive overview of the literature and the historical development of the problem, see  ~\cite{bottcher2013analysis,guille79,simon2005orthogonal}. Despite its relatively advanced age, the subject still attracts interest among mathematicians~\cite{guo2024first,MO24,nikolski2020szeg}.

The purpose of the present paper is both conceptual and concrete. Conceptually, we identify a common mechanism behind a broad class of first Szeg\H{o} limit theorems. We formulate an abstract theory for Toeplitz operators associated with normalized continuous Parseval frames. The hypotheses consist essentially of two properties: a normalization condition for the underlying measures and an approximate-identity condition expressing localization of the normalized frame vectors. Once these properties are verified, the proof of the Szeg\H{o} asymptotic is independent of the particular geometry of the underlying space and follows from the same argument based on Berezin--Lieb inequalities and convergence of Berezin transforms. Thus our approach separates the geometric input, which is needed to verify localization of the kernels, from the operator-theoretic argument producing the asymptotic trace formula.

\subsection{The abstract theorem}

The common mechanism behind all of our results is a Szeg\H{o} theorem for
Toeplitz operators associated with normalized continuous Parseval frames.
Let $X$ be a locally compact Hausdorff space and, for each parameter
$\alpha$, let $\mathcal{H}_\alpha$ be a Hilbert space equipped with a normalized
continuous Parseval frame
\[
  \{k_x^\alpha:x\in X\}\subset \mathcal{H}_\alpha
\]
and frame Radon measure $\nu_\alpha$.  Thus $\|k_x^\alpha\|_\alpha=1$ and
\[
  f=\int_X \langle f,k_x^\alpha\rangle_\alpha k_x^\alpha\,d\nu_\alpha(x),
  \qquad f\in \mathcal{H}_\alpha,
\]
where the integral is understood weakly.  We make two assumptions.
First, the frame measures have a common scale:
\begin{equation}
  d\nu_\alpha(x)=c(\alpha)\,d\nu(x),
  \label{intro:measure-scaling}
\end{equation}
where $\nu$ is a Radon measure independent of $\alpha$.  Second, the normalized frame
vectors become asymptotically localized: for every $x\in X$ and every
neighborhood $U$ of $x$,
\begin{equation}
  \int_{X\setminus U}
  \bigl|\langle k_x^\alpha,k_y^\alpha\rangle_\alpha\bigr|^2
  \,d\nu_\alpha(y)\longrightarrow 0
  \qquad (\alpha\to\infty).
  \label{intro:localization}
\end{equation}

For a real-valued symbol $\sigma\in L^1(X,\nu)$, set
\[
  T_\sigma^\alpha f
  =\int_X \sigma(x)\langle f,k_x^\alpha\rangle_\alpha
  k_x^\alpha\,d\nu_\alpha(x).
\]
The associated Berezin transform is
\[
  \widetilde\sigma_\alpha(x)
  =\int_X \sigma(y)
  \bigl|\langle k_x^\alpha,k_y^\alpha\rangle_\alpha\bigr|^2
  \,d\nu_\alpha(y).
\]

\medskip
\noindent\textbf{Main abstract theorem.}
Assumptions \eqref{intro:measure-scaling} and
\eqref{intro:localization} imply that
$\widetilde\sigma_\alpha\to\sigma$ in the relevant $L^p$ spaces.  Combined
with the Berezin--Lieb inequalities, this gives the following two trace
formulas.

\smallskip
\noindent\textbf{\emph{Compact form.}}
Suppose that $X$ is compact, that each $\mathcal{H}_\alpha$ is finite-dimensional, and
that $c(\alpha)=\dim \mathcal{H}_\alpha$.  If either
\begin{enumerate}
  \item $\sigma$ is bounded and $\psi$ is continuous on
  $[\operatorname*{ess\,inf}\sigma,\operatorname*{ess\,sup}\sigma]$, or
  \item $\sigma\geq 0$, $\psi\in C([0,\infty))$, and
  $\lim_{t\to\infty}\psi(t)/t$ exists and is finite,
\end{enumerate}
then
\begin{equation}
  \frac{1}{c(\alpha)}\operatorname{Tr}\bigl(\psi(T_\sigma^\alpha)\bigr)
  \longrightarrow
  \int_X\psi(\sigma(x))\,d\nu(x).
  \label{intro:compact-szego}
\end{equation}

\smallskip
\noindent\textbf{\emph{Locally compact form.}}
Suppose that
$\sigma,\eta\in L^1(X,\nu)\cap L^2(X,\nu)$ are real-valued,
$\eta\geq0$, and $T_\sigma^\alpha$ and $T_\eta^\alpha$ commute for every
$\alpha$.  If either
\begin{enumerate}
  \item $\sigma$ is bounded and $\psi$ is continuous on
  $[\operatorname*{ess\,inf}\sigma,\operatorname*{ess\,sup}\sigma]$, or
  \item $\sigma\geq0$, $\psi\in C([0,\infty))$, and
  $\lim_{t\to\infty}\psi(t)/t$ exists and is finite,
\end{enumerate}
then
\begin{equation}
  \frac{1}{c(\alpha)}
  \operatorname{Tr}\bigl(T_\eta^\alpha\psi(T_\sigma^\alpha)\bigr)
  \longrightarrow
  \int_X\eta(x)\psi(\sigma(x))\,d\nu(x).
  \label{intro:local-szego}
\end{equation}

The point is that all geometric information is confined to verifying
\eqref{intro:localization}.  Once this is done, the operator-theoretic proof
of \eqref{intro:compact-szego} and \eqref{intro:local-szego} is identical
across the examples.  We now record representative consequences for compact
Abelian groups, weighted Bergman spaces, and locally compact Abelian groups.

Concretely, the framework yields new Szeg\H{o} limit theorems as well as unified proofs and extensions of several known ones. In particular, we obtain a Szeg\H{o}-type limit theorem of this form for Toeplitz operators on weighted Bergman spaces of the unit ball. We also extend the compact Abelian-group setting of B\'edos~\cite{bedos1996folner} to the class of symbols considered below, recover the recent trace theorem of Guo, Li, and Zhou~\cite{guo2024first} as a special case, and obtain a corresponding theorem for locally compact Abelian groups. Further applications recover first-order asymptotics for Gabor localization and Paley--Wiener operators and yield eigenvalue-distribution consequences related to classical results of Landau--Widom and Lindholm.

We now state three representative consequences of the abstract theory.

\subsection{Weighted Bergman spaces}
 Our main application is a Szeg\H{o} limit theorem for Toeplitz operators on weighted Bergman spaces and other locally compact settings. 
The weighted Bergman spaces, denoted $A_\alpha^2(\mathbb{B}^{n})$ with $\alpha > -1$, consist of holomorphic functions on the unit ball $\mathbb{B}^{n}$ of $\mathbb{C}^n$ that are square-integrable with respect to the measure $d\mu_\alpha(z)=\frac{\Gamma(n+\alpha+1)}{\pi^{n}\Gamma(\alpha+1)}(1-|z|^2)^\alpha \, dV(z)$, 
where $dV(z)$ is the volume measure on $\mathbb{C}^n$. 
It is well known (see, e.g., \cite{Z90}) that $A_\alpha^2(\mathbb{B}^{n})$ is a reproducing kernel Hilbert space (RKHS), 
and that Toeplitz operators on these spaces take the form:

\begin{align*}
    T_{\sigma}^{\alpha}f=\frac{\Gamma(n+\alpha+1)}{\pi^{n}\Gamma(\alpha+1)}\int_{\mathbb{B}^{n}}\sigma(z)\ip{f}{k_{z}^{\alpha}}k_{z}^{\alpha}(1-|z|^{2})^{-n-1}dV(z), 
\end{align*}
where $k_{z}^{\alpha}(w)$ are the kernel sections. We prove the following analog of \eqref{Szegocont}. 


\begin{thm}\label{Sze2} Let $\sigma:\mathbb{B}^{n}\to \R$ be such that either
\begin{itemize}
    \item[(a)] $\sigma$ is bounded, $\sigma\in L^{1}(\mathbb{B}^{n},(1-|z|^{2})^{-n-1}dV(z))$, and $\psi:\R \to \R$ is continuous on $[\mathrm{ess\,inf}\,\sigma,\mathrm{ess\,sup}\,\sigma]$, or 
    \item[(b)] $\sigma\geq0$ and
\[
\sigma\in
L^1\!\left(B_n,(1-|z|^2)^{-n-1}dV(z)\right)
\cap
L^2\!\left(B_n,(1-|z|^2)^{-n-1}dV(z)\right),
\]
and $\psi\in C([0,\infty))$ satisfies
\[
\lim_{x\to\infty}\frac{\psi(x)}{x}\in\mathbb{R}.
\]
\end{itemize}
Then, 
    \begin{align*}
        \frac{\pi^{n}\Gamma(\alpha+1)}{\Gamma(n+\alpha+1)}\tr\left(T_{\sigma}^{\alpha}\psi\left(T_{\sigma}^{\alpha}\right)\right)\to\int_{\mathbb{B}^{n}}\sigma(z)\psi\left(\sigma(z)\right)(1-|z|^{2})^{-n-1}dV(z),
    \end{align*}
    as $\alpha\to\infty$.
\end{thm}

\noindent Our analog of \eqref{Szegocont} has slightly different form due to the non-compactness of the underlying domain $\mathbb{B}^{n}$. 

A corresponding \(n=1\) result was established in the first arXiv
version of the present work, posted in November 2024. Subsequently,
Ram{\'i}rez Monta{\~n}o~\cite{RamirezMontano} proved a related
Szeg\H{o} limit theorem on \(\mathbb B^n\) for smooth compactly supported
symbols and, more generally, for symbols supported on isotropic or
co-isotropic submanifolds. In the absolutely continuous case, his result
overlaps with the present theorem only for smooth compactly supported
symbols. By contrast, Theorem~\ref{Sze2} allows bounded \(L^1\) symbols, and
non-negative \(L^1\cap L^2\) symbols, without smoothness or compact-support
assumptions. The results are therefore complementary: the present theorem
extends the classical-symbol result in regularity, while
\cite{RamirezMontano} also treats singular symbol measures.

Besides providing a new result in this setting, the Bergman example illustrates the advantage of the abstract approach: the same proof applies despite the failure of the geometric assumptions underlying the standard Widom-type method.

\subsection{Compact Abelian groups}

\begin{thm}[Szeg\H{o} Limit on Compact Abelian Groups]
\label{thm:compact_ab_group}
Let $G$ be a compact metrizable Abelian group with normalized Haar measure $\mu$ and let $\{\Gamma_N\}$ be a sequence of finite subsets of $\widehat G$ satisfying the F\o lner condition
\[
\frac{|(\gamma\Gamma_N)\triangle\Gamma_N|}{|\Gamma_N|}
\longrightarrow0
\qquad \gamma\in\widehat G.
\]
Let $\mathcal{K}_{N}=\text{span}\{\Gamma_{N}\}$ and $\sigma \in L^1(G)$ be a real-valued, non-negative function. Let $T_\sigma^N : \mathcal{K}_N \to \mathcal{K}_N$ denote the corresponding Toeplitz operator
\[
T_\sigma^N f = P_N (\sigma f), \quad f \in \mathcal{K}_N,\]
where $P_N : L^2(G) \to \mathcal{K}_N$ is the orthogonal projection. Then, for any continuous function $\psi : [0,\infty) \to \mathbb{R}$ such that $\psi(x)/x$ has finite limit as $x \to \infty$, we have
\[
\lim_{N\to\infty} \frac{1}{|\Gamma_N|} \mathrm{Tr} \, \psi(T_\sigma^N) = \int_G \psi(\sigma(x)) \, d\mu(x).
\]
\end{thm}

An important feature of our approach is that it also yields the
determinant formulation of the first Szeg\H{o} limit theorem.  Although
the logarithm is not covered directly by our continuous functional
calculus results because of its singularity at the origin, the
determinant asymptotic follows from the same abstract framework by a
regularization argument.  More precisely, for a non-negative
integrable symbol $\sigma$ in the compact setting, we obtain
\[ \lim_{N\to\infty} \frac{1}{|\Gamma_N|} \log\det T_\sigma^N = \int_G \log\sigma(x)\,d\mu(x), \]
where the right-hand side, and hence the limit, is allowed to take the
value $-\infty$.  Thus no assumption that the symbol be bounded away
from zero is required.  In particular, our abstract framework yields
both the trace and determinant formulations of the first Szeg\H{o}
limit theorem.

\subsection{Locally compact Abelian groups}

The same framework also produces a non-compact analogue of the compact-group theorem.

\begin{thm}[Szeg\H{o} theorem for LCA groups]\label{thm:szego_lca}
Let $G$ be a locally compact Abelian group with Haar measure $m$
and dual group $\widehat{G}$ with Haar measure $\widehat{m}$.
Let $\{F_N\}_{N\ge1}\subset \widehat{G}$ be a F\o lner sequence and define
\[
\mathcal{H}_N
=
\{f\in L^2(G): \operatorname{supp}(\widehat{f})\subset F_N\}.
\]
Let $c(N)=\widehat{m}(F_N)$.

For $\sigma\in L^1(G,m)$ real-valued,
let $T_\sigma^N$ be the Toeplitz operator
\[
T_\sigma^N f = P_N(\sigma f).
\]

If either

\begin{itemize}
\item[(a)] $\sigma$ is bounded, $\sigma\in L^{1}(G)$, and
$\psi$ is continuous on
$[\mathrm{ess\,inf}\,\sigma,\mathrm{ess\,sup}\,\sigma]$, or
\item[(b)] $\sigma\ge0$, $\sigma\in L^{2}(G)$, and
$\psi\in C([0,\infty))$ such that the limit
$\lim_{x\to\infty}\psi(x)/x$ exists,
\end{itemize}

then
\[
\lim_{N\to\infty}
\frac{1}{c(N)}
\operatorname{tr}\!\big(
T_\sigma^N \psi(T_\sigma^N)
\big)
=
\int_G \sigma(x)\psi(\sigma(x))\, dm(x).
\]
\end{thm}

This may be regarded as a locally compact counterpart of the compact Abelian-group Szeg\H{o} theorem. More importantly for our purposes, both results arise from exactly the same abstract argument.

\subsection{Relation to previous work and scope}

The results above are first-order trace formulae.  Their purpose is to isolate a localization principle that produces the leading term in several different settings; they are not intended to subsume the finer, model-dependent asymptotics available in all of those settings.

For compact abelian groups, B\'edos~\cite{bedos1996folner} proved a Szeg\H{o}-type trace formula for bounded real-valued symbols along F{\o}lner sequences. For \(G=\mathbb T^d\), \(1\leq d\leq\infty\), Guo, Li, and Zhou~\cite{guo2024first} extended the trace formula to nonnegative \(L^1\)-symbols and continuous test functions of at most linear growth (motivated in part by questions of Nikolski and Pushnitski~\cite{nikolski2020szeg}); they also obtained determinant formulae and studied certain non-F{\o}lner
truncations. Our Theorem~\ref{thm:compact_ab_group} and Corollary~\ref{DetGroup} recover the trace and determinant part of their results in the torus case and extend the underlying compact group to an arbitrary compact
metrizable abelian group.  It does not address their 
non-F{\o}lner results, so those results are complementary to ours.

The comparison with the weighted Bergman theorem is of a different nature.
Ramirez Monta{\~n}o~\cite{RamirezMontano} treats positive smooth,
compactly supported symbols associated with isotropic and co-isotropic
submanifolds; in particular, his setting includes geometrically defined
singular measure symbols.  By contrast, Theorem~\ref{Sze2} concerns absolutely
continuous symbols and permits the rough symbol classes stated there, with
no smoothness or compact-support assumption.  Thus, in the common case of a
positive smooth compactly supported density, the two results give compatible
first-order formulae after the normalizations are matched.  Neither result
contains the other: the former allows singular geometric symbols, whereas
the latter allows rough diffuse symbols.

In the Gabor setting, the theorem of Feichtinger and
Nowak~\cite{FN01} assumes \(b\in L^1\) with \(0\leq b\leq1\).
Since this implies \(b\in L^1\cap L^2\), Corollary~\ref{FeiNow} recovers that
first-order asymptotic and permits the broader nonnegative symbol class
covered by our local theorem, together with the corresponding
linear-growth class of test functions.  Likewise, the Paley--Wiener
specialization recovers the leading Landau--Widom counting law
\cite{landau1980eigenvalue}.  It does not recover the logarithmic second-order term,
which depends on boundary information beyond the localization criterion used
here.

\section{Berezin--Lieb type inequalities}

We record two Berezin--Lieb type inequalities that will be used later. Let $\mathfrak{S}_{1}(\mathcal{H}_{\alpha})$ denote the trace-class operators.  

\subsection*{The classical inequality}

\begin{prop}\label{BL1}
Let $\sigma\in L^1(X,\nu_\alpha)$ be real-valued,
and let $\psi:\mathbb{R}\to\mathbb{R}$ be convex. Assume that
\[
\psi(T_\sigma^\alpha)\in\mathfrak{S}_1(\mathcal H_\alpha)
\qquad\text{and}\qquad
\psi(\sigma)\in L^1(X,\nu_\alpha).
\]
Then
\begin{equation}\label{BL-classical}
\int_X
\psi\bigl(\widetilde{\sigma}^{\,\alpha}(x)\bigr)
\,d\nu_\alpha(x)
\leq
\operatorname{Tr}\bigl(\psi(T_\sigma^\alpha)\bigr)
\leq
\int_X\psi\bigl(\sigma(x)\bigr)\,d\nu_\alpha(x).
\end{equation}
\end{prop}

\begin{proof} Since $\sigma\in L^1(X,\nu_\alpha)$, the operator
$T_\sigma^\alpha$ is trace class and therefore compact. Since
$\sigma$ is real-valued, it is self-adjoint. Hence there exists an
orthonormal eigenbasis $(\phi_j)_j$ of $\mathcal H_\alpha$, with
corresponding real eigenvalues $(\lambda_j)_j$.

Since the frame is normalized and Parseval,
\[
\sum_j
\bigl|\langle \phi_j,k_x^\alpha\rangle_\alpha\bigr|^2=1
\]
for every $x\in X$. Moreover,
\[
\widetilde{\sigma}^{\,\alpha}(x)
=
\left\langle T_\sigma^\alpha k_x^\alpha,k_x^\alpha
\right\rangle_\alpha
=
\sum_j \lambda_j
\bigl|\langle \phi_j,k_x^\alpha\rangle_\alpha\bigr|^2.
\]
Jensen's inequality therefore gives
\[
\psi\bigl(\widetilde{\sigma}^{\,\alpha}(x)\bigr)
\leq
\sum_j \psi(\lambda_j)
\bigl|\langle \phi_j,k_x^\alpha\rangle_\alpha\bigr|^2.
\]
Integrating and using the Parseval frame identity,
\[
\int_X
\psi\bigl(\widetilde{\sigma}^{\,\alpha}(x)\bigr)
\,d\nu_\alpha(x)
\leq
\sum_j\psi(\lambda_j)
=
\operatorname{Tr}\bigl(\psi(T_\sigma^\alpha)\bigr).
\]

For the reverse inequality, we have
\[
\lambda_j
=
\left\langle T_\sigma^\alpha\phi_j,\phi_j\right\rangle_\alpha
=
\int_X \sigma(x)
\bigl|\langle \phi_j,k_x^\alpha\rangle_\alpha\bigr|^2
\,d\nu_\alpha(x).
\]
Applying Jensen's inequality once more,
\[
\psi(\lambda_j)
\leq
\int_X \psi(\sigma(x))
\bigl|\langle \phi_j,k_x^\alpha\rangle_\alpha\bigr|^2
\,d\nu_\alpha(x).
\]
Summing over $j$ and using Parseval again yields
\[
\operatorname{Tr}\bigl(\psi(T_\sigma^\alpha)\bigr)
\leq
\int_X \psi(\sigma(x))\,d\nu_\alpha(x).
\]
This proves the result.
\end{proof}

Under Assumption~(I), dividing \eqref{BL-classical} by $c(\alpha)$ gives

\begin{equation}\label{renorm1}
\int_X
\psi\bigl(\widetilde{\sigma}^{\,\alpha}(x)\bigr)\,d\nu(x)
\leq
\frac{1}{c(\alpha)}
\operatorname{Tr}\bigl(\psi(T_\sigma^\alpha)\bigr)
\leq
\int_X\psi\bigl(\sigma(x)\bigr)\,d\nu(x).
\end{equation}

\subsection*{A weighted inequality}

\begin{prop}\label{BLweighted}
Let $\sigma,\eta\in L^1(X,\nu_\alpha)$, with
$\sigma$ real-valued and $\eta\geq 0$. Assume that
$T_\sigma^\alpha$ and $T_\eta^\alpha$ commute. Let
$\psi:\mathbb{R}\to\mathbb{R}$ be convex. Then

\begin{equation}\label{BL-weighted}
\int_X
\eta(x)\,
\psi\bigl(\widetilde{\sigma}^{\,\alpha}(x)\bigr)
\,d\nu_\alpha(x)
\leq
\operatorname{Tr}\!\left(
T_\eta^\alpha\psi(T_\sigma^\alpha)
\right)
\leq
\int_X
\widetilde{\eta}^{\,\alpha}(x)\,
\psi\bigl(\sigma(x)\bigr)
\,d\nu_\alpha(x),
\end{equation}
whenever the right-hand side is finite.
\end{prop}

\begin{proof}
Since $\sigma,\eta\in L^1(X,\nu_\alpha)$, the operators
$T_\sigma^\alpha$ and $T_\eta^\alpha$ are trace class. Moreover,
$T_\sigma^\alpha$ is self-adjoint and $T_\eta^\alpha$ is positive.
Since they commute, they admit a common orthonormal eigenbasis
$(\phi_j)_j$, with
\[
T_\sigma^\alpha\phi_j=\lambda_j\phi_j,
\qquad
T_\eta^\alpha\phi_j=\mu_j\phi_j,
\qquad
\mu_j\geq 0.
\]
Thus
\[
\operatorname{Tr}\!\left(
T_\eta^\alpha\psi(T_\sigma^\alpha)
\right)
=
\sum_j\mu_j\psi(\lambda_j).
\]

For each $x\in X$, Parseval's identity gives
\[
\widetilde{\sigma}^{\,\alpha}(x)
=
\sum_j\lambda_j
\bigl|\langle\phi_j,k_x^\alpha\rangle_\alpha\bigr|^2,
\]
and
\[
\sum_j
\bigl|\langle\phi_j,k_x^\alpha\rangle_\alpha\bigr|^2
=1.
\]
Therefore, Jensen's inequality implies
\[
\psi\bigl(\widetilde{\sigma}^{\,\alpha}(x)\bigr)
\leq
\sum_j\psi(\lambda_j)
\bigl|\langle\phi_j,k_x^\alpha\rangle_\alpha\bigr|^2.
\]
Multiplying by $\eta(x)$ and integrating gives
\[
\begin{aligned}
\int_X\eta(x)
\psi\bigl(\widetilde{\sigma}^{\,\alpha}(x)\bigr)
\,d\nu_\alpha(x)
&\leq
\sum_j\psi(\lambda_j)
\int_X\eta(x)
\bigl|\langle\phi_j,k_x^\alpha\rangle_\alpha\bigr|^2
\,d\nu_\alpha(x)\\
&=
\sum_j\mu_j\psi(\lambda_j).
\end{aligned}
\]

For the reverse inequality, we have
\[
\lambda_j
=
\int_X\sigma(x)
\bigl|\langle\phi_j,k_x^\alpha\rangle_\alpha\bigr|^2
\,d\nu_\alpha(x).
\]
Another application of Jensen's inequality yields
\[
\psi(\lambda_j)
\leq
\int_X\psi(\sigma(x))
\bigl|\langle\phi_j,k_x^\alpha\rangle_\alpha\bigr|^2
\,d\nu_\alpha(x).
\]
Multiplying by $\mu_j$ and summing over $j$, we obtain
\[
\begin{aligned}
\sum_j\mu_j\psi(\lambda_j)
&\leq
\int_X\psi(\sigma(x))
\sum_j\mu_j
\bigl|\langle\phi_j,k_x^\alpha\rangle_\alpha\bigr|^2
\,d\nu_\alpha(x)\\
&=
\int_X
\widetilde{\eta}^{\,\alpha}(x)
\psi(\sigma(x))\,d\nu_\alpha(x).
\end{aligned}
\]
This proves the result.
\end{proof}

Dividing by $c(\alpha)$ gives
\begin{equation}\label{renorm2}
\int_X
\eta(x)\,
\psi\bigl(\widetilde{\sigma}^{\,\alpha}(x)\bigr)
\,d\nu(x)
\leq
\frac{1}{c(\alpha)}
\operatorname{Tr}\!\left(
T_\eta^\alpha\psi(T_\sigma^\alpha)
\right)
\leq
\int_X
\widetilde{\eta}^{\,\alpha}(x)\,
\psi\bigl(\sigma(x)\bigr)
\,d\nu(x).
\end{equation}

\medskip

If $\psi$ is concave, the inequalities in
\eqref{BL-classical} and \eqref{BL-weighted} reverse.

\medskip

\noindent
Both inequalities are particularly useful after normalization.
If $\widetilde{\sigma}^\alpha$ converges to $\sigma$
in a suitable sense and the family
$\{\widetilde{\eta}^\alpha\}$ admits an analogous limit,
then \eqref{renorm1} and \eqref{renorm2}
yield Szeg\H{o}-type trace asymptotics by a standard squeezing argument
combined with an appropriate convergence theorem.
\section{Convergence of the Berezin transform}

The goal of this section is to show that an approximate identity
condition implies convergence in $L^p$ (hence in measure)
of the Berezin transform, which is the key input for the
Szeg\H{o} limit theorems of Sections~4 and~5.

\medskip

\noindent
\begin{itemize}
    \item \textbf{Assumption (II)} For every $x\in X$ and every neighborhood $N_x$ of $x$,
\[
    \int_{N_x^c}
    |\langle k_x^\alpha,k_y^\alpha\rangle|^2
    \, d\nu_\alpha(y)
    \;\longrightarrow\; 0
    \qquad (\alpha\to\infty).
\]
\end{itemize}

\begin{lem}
For $1 \le p \le \infty$, the Berezin transform
$\sigma \mapsto \widetilde{\sigma}^\alpha$
defines a linear contraction on $L^p(X,\nu)$,
uniformly in $\alpha$.
\end{lem}

\begin{proof}
The $L^\infty$ bound follows immediately from
\[
    |\widetilde{\sigma}^\alpha(x)|
    \le
    \|\sigma\|_\infty
    \int_X
    |\langle k_x^\alpha,k_y^\alpha\rangle|^2
    d\nu_\alpha(y)
    =
    \|\sigma\|_\infty.
\]

For $p=1$, Fubini's theorem and Assumption (I) give
\[
    \|\widetilde{\sigma}^\alpha\|_1
    \le
    \int_X |\sigma(y)|\, d\nu(y)
    =
    \|\sigma\|_1.
\]

The general case follows by interpolation.
\end{proof}

\begin{lem}
If $\sigma \in C_c(X)$, then
$\widetilde{\sigma}^\alpha(x) \to \sigma(x)$
for every $x\in X$.
\end{lem}

\begin{proof}
Fix $x\in X$ and $\varepsilon>0$.
By continuity, there exists a neighborhood $N_x$
such that $|\sigma(y)-\sigma(x)|<\varepsilon$
for $y\in N_x$.

Write
\[
    \widetilde{\sigma}^\alpha(x)-\sigma(x)
    =
    \int_X
    (\sigma(y)-\sigma(x))
    |\langle k_x^\alpha,k_y^\alpha\rangle|^2
    d\nu_\alpha(y).
\]
Splitting the integral over $N_x$ and $N_x^c$,
the contribution on $N_x$ is bounded by $\varepsilon$,
while the contribution on $N_x^c$
tends to $0$ by Assumption~(II)
since $\sigma$ is bounded.
\end{proof}

\begin{prop}\label{Lp}
Let $1 \le p < \infty$ and $\sigma \in L^p(X,\nu)$.
Then
\[
    \widetilde{\sigma}^\alpha \to \sigma
    \quad \text{in } L^p(X,\nu).
\]
In particular, $\widetilde{\sigma}^\alpha \to \sigma$
in measure.
\end{prop}

\begin{proof}
First assume $\sigma \in C_c(X)$.
By the previous lemma,
$\widetilde{\sigma}^\alpha(x) \to \sigma(x)$ pointwise.
Moreover, Jensen’s inequality gives
\[
    |\widetilde{\sigma}^\alpha(x)|^p
    \le
    \widetilde{|\sigma|^p}^{\,\alpha}(x).
\]
Hence
\[
    |\widetilde{\sigma}^\alpha(x)-\sigma(x)|^p
    \le
    2^{p-1}
    \big(
        \widetilde{|\sigma|^p}^{\,\alpha}(x)
        + |\sigma(x)|^p
    \big).
\]
The right-hand side has constant integral,
since
\[
    \int_X
    \widetilde{|\sigma|^p}^{\,\alpha}(x)
    d\nu(x)
    =
    \int_X |\sigma(x)|^p d\nu(x).
\]
Since $|\sigma|^p\in C_c(X)$, Lemma~3.2 also gives
\[
\widetilde{|\sigma|^p}^{\,\alpha}(x)\longrightarrow |\sigma(x)|^p
\qquad\text{for every }x\in X.
\]
Thus the generalized dominated convergence theorem applies to
\[
f_\alpha:=\bigl|\widetilde{\sigma}^{\alpha}-\sigma\bigr|^p
\quad\text{and}\quad
g_\alpha:=2^{p-1}\left(
\widetilde{|\sigma|^p}^{\,\alpha}+|\sigma|^p
\right).
\]
Indeed, $0\leq f_\alpha\leq g_\alpha$, $f_\alpha\to0$ pointwise,
$g_\alpha\to2^p|\sigma|^p$ pointwise, and the preceding identity gives
\[
\int_X g_\alpha\,d\nu
=
\int_X 2^p|\sigma|^p\,d\nu.
\]
Consequently,
\[
\bigl\|\widetilde{\sigma}^{\alpha}-\sigma\bigr\|_p\longrightarrow0.
\]

For general $\sigma \in L^p$,
approximate by $\psi \in C_c(X)$ and use
the uniform $L^p$ boundedness of the Berezin transform.
\end{proof}

The preceding results apply to Toeplitz operators with integrable symbols, which may be viewed as Toeplitz operators whose symbol is an absolutely continuous measure on phase space. More generally, one may define Toeplitz operators whose symbol is a (finite) complex Borel measure.

\begin{definition}
Let $\mu:\mathcal{B}(X)\to\C$ be a complex-valued Borel measure. The (weakly defined) Toeplitz operator with symbol $\mu$ is given by
\[
T_\mu f
= \int_X \langle f,k_x^\alpha\rangle_\alpha \, k_x^\alpha \, d\mu(x).
\]
\end{definition}

Although the Szeg\H{o}-type limit theorems do not extend in full generality to this setting, the associated Berezin transforms remain well defined. However, the separability of the measure $\nu_{\alpha}$ no longer makes sense in this circumstance. Instead, we replace Assumptions (I) and (II) with the following assumption, which still capture this essence:
\begin{itemize}
    \item \textbf{Assumption (III)} For each $y\neq x$, $|\ip{k_{x}^{\alpha}}{k_{y}^{\alpha}}_{\alpha}|^{2}\to 0$ as $\alpha\to\infty$. 
\end{itemize}
From this assumption, we obtain the following. 

\begin{prop}
Let $\mu$ be a complex Borel measure with $|\mu|(X)<\infty$, and let $T_\mu:\mathcal{H}_\alpha\to\mathcal{H}_\alpha$ denote the corresponding Toeplitz operators. Then, for every $x\in X$,
\[
\widetilde{\mu}^\alpha(x)
:= \langle T_\mu k_x^\alpha, k_x^\alpha\rangle_\alpha
\longrightarrow \mu(\{x\}),
\qquad \text{as } \alpha\to\infty.
\]
\end{prop}

\begin{proof}


Let $x\in X$. Define the measure $\hat{\mu}:=|\mu-\mu(\{x\})\delta_{x}|$. Then $\hat{\mu}$ is a finite measure. Proceeding along, let $U\subset X$ be an open neighborhood of $x$. Then, 
    \begin{align*}
        \left|\Tilde{\mu}^{\alpha}(x)-\mu(\{x\})\right|=&\left|\int_{X}\left|\ip{k_{x}^{\alpha}}{k_{y}^{\alpha}}_{\alpha}\right|^{2}d\mu(y)-\mu(\{x\})\int_{X}\left|\ip{k_{x}^{\alpha}}{k_{y}^{\alpha}}_{\alpha}\right|^{2}d\delta_{x}(y)\right|\\
        \leq&\int_{U}\left|\ip{k_{x}^{\alpha}}{k_{y}^{\alpha}}_{\alpha}\right|^{2}d\hat{\mu}(y)+\int_{U^{c}}\left|\ip{k_{x}^{\alpha}}{k_{y}^{\alpha}}_{\alpha}\right|^{2}d\hat{\mu}(y).
    \end{align*}
    Let $\varepsilon>0$. Now, as $\left|\ip{k_{x}^{\alpha}}{k_{y}^{\alpha}}_{\alpha}\right|^{2}\leq 1$, and $\hat{\mu}$ is a finite measure with $\hat{\mu}(\{x\})=0$, we may choose $U$ sufficiently small such that $\hat{\mu}(U)<\varepsilon/2$. Thus, 
    \begin{align*}
        \int_{U}\left|\ip{k_{x}^{\alpha}}{k_{y}^{\alpha}}_{\alpha}\right|^{2}d\hat{\mu}(y)\leq\hat{\mu}(U)<\varepsilon/2.
    \end{align*}
    Likewise, $\left|\ip{k_{x}^{\alpha}}{k_{y}^{\alpha}}_{\alpha}\right|^{2}\leq 1$ on $U^{c}$, and since $\hat{\mu}(\{x\})=0$ and Assumption (III), we have that $\left|\ip{k_{x}^{\alpha}}{k_{y}^{\alpha}}_{\alpha}\right|^{2}\to 0$ $\hat{\mu}$-a.e. as $\alpha\to\infty$. Thus, by the Bounded Convergence Theorem, there is an $\alpha$ sufficiently large such that 
    \begin{align*}
        \int_{U^{c}}\left|\ip{k_{x}^{\alpha}}{k_{y}^{\alpha}}_{\alpha}\right|^{2}d\hat{\mu}(y)<\varepsilon/2.
    \end{align*}
    Combining the results gives the desired limit. 
\end{proof}

\begin{prop}\label{lkc}
Let $\mu$ be a finite complex Radon measure on $X$. Under
Assumption~(III), for every $x\in X$,
\[
\int_X
\langle k_x^\alpha,k_y^\alpha\rangle_\alpha\,d\mu(y)
\longrightarrow
\mu(\{x\})
\qquad
(\alpha\to\infty).
\]
\end{prop}

\begin{proof}
Set
\[
A_\alpha(x)
=
\int_X
\langle k_x^\alpha,k_y^\alpha\rangle_\alpha\,d\mu(y)
\]
and
\[
\widehat\mu
=
\left|\mu-\mu(\{x\})\delta_x\right|.
\]
Since $\langle k_x^\alpha,k_x^\alpha\rangle_\alpha=1$, we have
\[
A_\alpha(x)-\mu(\{x\})
=
\int_X
\langle k_x^\alpha,k_y^\alpha\rangle_\alpha
\,d\bigl(\mu-\mu(\{x\})\delta_x\bigr)(y).
\]

Let $\varepsilon>0$. Choose a neighborhood $U$ of $x$ such that
\[
\widehat\mu(U)<\frac{\varepsilon}{2}.
\]
Since the kernels are normalized,
\[
\left|\langle k_x^\alpha,k_y^\alpha\rangle_\alpha\right|
\leq 1.
\]
Therefore,
\[
\begin{aligned}
|A_\alpha(x)-\mu(\{x\})|
&\leq
\int_U
\left|\langle k_x^\alpha,k_y^\alpha\rangle_\alpha\right|
\,d\widehat\mu(y) \\
&\quad+
\int_{U^c}
\left|\langle k_x^\alpha,k_y^\alpha\rangle_\alpha\right|
\,d\widehat\mu(y).
\end{aligned}
\]
The first term is smaller than $\varepsilon/2$. For $y\neq x$,
Assumption~(III) implies
\[
\left|\langle k_x^\alpha,k_y^\alpha\rangle_\alpha\right|
\longrightarrow 0.
\]
Hence the second term tends to zero by dominated convergence. Thus
\[
\limsup_{\alpha\to\infty}
|A_\alpha(x)-\mu(\{x\})|
\leq \frac{\varepsilon}{2}.
\]
Since $\varepsilon$ is arbitrary, the result follows.
\end{proof}

\section{Abstract Szeg\H{o} Limit Theorem for Compact Spaces}

Throughout this section, let $\sigma \in L^1(X,\nu)$ be real-valued, and $\nu$ is a finite measure. We also assume that all the Hilbert spaces $H_\alpha$ are finite-dimensional. Under these assumptions the associated Toeplitz operators $T_\sigma^\alpha$
are trace-class, self-adjoint and satisfy
\[
\widetilde{\sigma}^\alpha \longrightarrow \sigma
\quad\text{in $\nu$-measure as } \alpha\to\infty.
\]

We begin with the convex case, which forms the core of the argument.

\begin{thm}\label{thm:convex_szego}
\leavevmode
\begin{itemize}
\item[(i)] Suppose $\sigma \ge 0$ and $\psi \in C([0,\infty))$
is convex such that the limit 
\[
\lim_{x\to\infty} \frac{\psi(x)}{x}
\]
exists. Then
\[
\lim_{\alpha\to\infty}
\frac{1}{c(\alpha)}\tr\!\big(\psi(T_\sigma^\alpha)\big)
=
\int_X \psi(\sigma(x))\,d\nu(x).
\]

\item[(ii)] If $\sigma \in L^\infty(X,\nu)$ and
$\psi$ is convex and continuous on
$[\text{ess inf } \sigma, \text{ess sup }\sigma]$, then the same limit holds.
\end{itemize}
\end{thm}

\begin{proof}
We prove (i); part (ii) follows by a simpler argument.

For $\psi(x)=x$ the result holds since
\[
\frac{1}{c(\alpha)}\tr(T_\sigma^\alpha)
=
\int_X \sigma\, d\nu,
\]
which is independent of $\alpha$. 

Write
\[
\psi(x)=\beta x + \psi_0(x),
\qquad
\text{where } \frac{\psi_0(x)}{x}\to 0.
\]
By linearity it suffices to treat the case
$\psi(x)/x\to 0$.

For convex $\psi$, Proposition~\ref{BL1} gives
\[
\int_X \psi(\widetilde{\sigma}^\alpha)\,d\nu_\alpha
\le
\tr(\psi(T_\sigma^\alpha))
\le
\int_X \psi(\sigma)\,d\nu_\alpha.
\]
Dividing by $c(\alpha)$ and using Assumption~(I),
\[
\int_X \psi(\widetilde{\sigma}^\alpha)\,d\nu
\le
\frac{1}{c(\alpha)}\tr(\psi(T_\sigma^\alpha))
\le
\int_X \psi(\sigma)\,d\nu.
\]

Since $\psi(x)/x\to 0$, for every $\varepsilon>0$
there exists $C_\varepsilon$ such that
\[
|\psi(x)| \le \varepsilon x + C_\varepsilon.
\]
By Proposition~\ref{Lp}, $\widetilde{\sigma}^{\alpha}\to\sigma$ in
$L^1(X,\nu)$, and hence the family
$\{\widetilde{\sigma}^{\alpha}\}_{\alpha}$ is uniformly integrable.
Since $\nu(X)<\infty$, the preceding estimate implies that the family
$\{\psi(\widetilde{\sigma}^{\alpha})\}_{\alpha}$ is uniformly
integrable. Moreover,
\[
\psi(\widetilde{\sigma}^{\alpha})\longrightarrow\psi(\sigma)
\]
in measure. Therefore, Vitali's convergence theorem gives
\[
\int_X\psi(\widetilde{\sigma}^{\alpha})\,d\nu
\longrightarrow
\int_X\psi(\sigma)\,d\nu .
\]

\medskip
\noindent
\textbf{Proof of (ii).}
If $\sigma\in L^\infty$, its essential range
is compact. A convex continuous function on a compact
interval is uniformly continuous.
The same Berezin--Lieb argument as above,
now without growth considerations,
yields the limit directly.
\end{proof}

\subsection{From Convex Functions to General Continuous Symbols}

We extend the Szeg\H{o} theorem from convex functions to general continuous functions via a decomposition into differences of convex functions.

\begin{lem}\label{lem:diff_convex}
Let $\psi$ be real-valued. 
\begin{itemize}
    \item[(i)] If $\psi\in C[a,b]$, then for every $\varepsilon>0$ there exist convex functions $\psi_1,\psi_2$ on $[a,b]$ such that 
    \[
        \|\psi-(\psi_1-\psi_2)\|_{\infty} < \varepsilon.
    \]
    \item[(ii)] If $\psi\in C_0[0,\infty)$, then for every $\varepsilon>0$ there exist convex functions $\psi_1,\psi_2\in C_0[0,\infty)$ such that
    \[
        \|\psi-(\psi_1-\psi_2)\|_\infty < \varepsilon.
    \]
\end{itemize}
\end{lem}

\begin{proof}
(i) For a polynomial $p$ on $[a,b]$, let $M = \inf_{x\in[a,b]} p''(x)$. If $M\geq 0$ then $p$ is convex itself, so we can trivially decompose $p=p-0$. If $M<0$ then 
\[
    p(x) = (p(x)-\frac{M}{2}x^2) - (-\frac{M}{2})x^2,
\]
where both terms are convex. Any $\psi\in C[a,b]$ can be uniformly approximated by polynomials; applying this decomposition yields the claim.

(ii) Let
\[
\mathcal C
:=
\{\varphi\in C_0([0,\infty)):
\varphi\text{ is convex}\}
\]
and
\[
\mathcal A
:=
\{\varphi_1-\varphi_2:
\varphi_1,\varphi_2\in\mathcal C\}.
\]

We first show that every $\varphi\in\mathcal C$ is non-negative and
non-increasing. If $0\leq s<t<u$, convexity gives
\[
\frac{\varphi(t)-\varphi(s)}{t-s}
\leq
\frac{\varphi(u)-\varphi(t)}{u-t}.
\]
Letting $u\to\infty$ and using $\varphi(u)\to 0$, we obtain
\[
\frac{\varphi(t)-\varphi(s)}{t-s}\leq 0.
\]
Thus $\varphi$ is non-increasing. Since $\varphi(t)\to 0$ as
$t\to\infty$, it follows that $\varphi\geq 0$.

We next show that $\mathcal C$ is closed under multiplication. Let
$f,g\in\mathcal C$. Then $fg\in C_0([0,\infty))$. Suppose that
$x\leq y$ and let $z=\lambda x+(1-\lambda)y$, where
$0\leq\lambda\leq 1$. Set
\[
a=f(x),\quad b=f(y),\quad c=g(x),\quad d=g(y).
\]
Since $f$ and $g$ are non-negative and non-increasing,
$a\geq b\geq 0$ and $c\geq d\geq 0$. By convexity,
\[
f(z)\leq \lambda a+(1-\lambda)b,
\qquad
g(z)\leq \lambda c+(1-\lambda)d.
\]
Moreover,
\[
\begin{aligned}
&\lambda ac+(1-\lambda)bd
\\
&\quad-
\bigl(\lambda a+(1-\lambda)b\bigr)
\bigl(\lambda c+(1-\lambda)d\bigr)
\\
&=
\lambda(1-\lambda)(a-b)(c-d)\geq 0.
\end{aligned}
\]
Consequently,
\[
f(z)g(z)
\leq
\lambda f(x)g(x)+(1-\lambda)f(y)g(y),
\]
so $fg$ is convex. Therefore $\mathcal C$ is closed under
multiplication.

It follows that $\mathcal A$ is an algebra. Indeed, if
$f=f_1-f_2$ and $g=g_1-g_2$, with $f_j,g_j\in\mathcal C$, then
\[
fg
=
(f_1g_1+f_2g_2)-(f_1g_2+f_2g_1)\in\mathcal A.
\]

The function
\[
\rho(t)=\frac{1}{1+t}
\]
belongs to $\mathcal C\subset\mathcal A$, is strictly decreasing, and
is nonzero at every point. Hence $\mathcal A$ separates points and
vanishes nowhere.

By the Stone--Weierstrass theorem for $C_0([0,\infty))$, $\mathcal A$
is uniformly dense in $C_0([0,\infty))$. Therefore, for every
$\varepsilon>0$, there exist convex functions
$\psi_1,\psi_2\in C_0([0,\infty))$ such that
\[
\|\psi-(\psi_1-\psi_2)\|_\infty<\varepsilon.
\]
This proves (ii).
\end{proof}

\subsection{General Szeg\H{o} Theorem for compact $X$}

\begin{thm}\label{thm:general_szego}
Let $\sigma\in L^1(X,\nu)$ be real-valued. Suppose $c(\alpha)=\dim \mathcal{H}_\alpha$ for all $\alpha$.
\begin{itemize}
    \item[(i)] If $\sigma$ is bounded, let $\psi$ be continuous on $[\text{ess inf } \sigma, \text{ess sup }\sigma]$.
    \item[(ii)] If $\sigma$ is unbounded and non-negative, let $\psi\in C[0,\infty)$ such that the limit $\lim_{x\to\infty} \psi(x)/x$ exists.
\end{itemize}

Then
\[
    \lim_{\alpha\to\infty} \frac{1}{c(\alpha)} \tr\big(\psi(T_\sigma^\alpha)\big)
    = \int_X \psi(\sigma(x)) \, d\nu(x).
\]
\end{thm}

\begin{proof}
\textbf{Case (i):} By Lemma~\ref{lem:diff_convex} (i), for every $\varepsilon>0$ there exist convex functions $\psi_1,\psi_2$ on $[\text{ess inf } \sigma, \text{ess sup }\sigma]$ such that
\[
    \|\psi-(\psi_1-\psi_2)\|_\infty < \varepsilon.
\]
Then
\[
    \frac{1}{c(\alpha)} \tr(\psi(T_\sigma^\alpha)) 
    = \frac{1}{c(\alpha)} \tr(\psi_1(T_\sigma^\alpha)) 
      - \frac{1}{c(\alpha)} \tr(\psi_2(T_\sigma^\alpha)) 
      + R_\alpha,
\]
where $|R_\alpha| \le \varepsilon$. Applying the convex Szeg\H{o} theorem to $\psi_1,\psi_2$ and letting $\varepsilon\to 0$ gives the result.

\textbf{Case (ii):} 
Let
\[
\beta:=\lim_{x\to\infty}\frac{\psi(x)}{x},
\qquad
r(x):=\psi(x)-\beta x.
\]
Then
\[
\frac{r(x)}{x}\longrightarrow 0
\qquad (x\to\infty).
\]
It is therefore enough to prove the result for $r$.

Fix $\varepsilon>0$. Choose $R>0$ such that
\[
|r(t)|\leq \varepsilon t
\qquad\text{for all }t\geq R.
\]
Let $\chi_R\in C_c([0,\infty))$ satisfy
\[
0\leq\chi_R\leq 1,\qquad
\chi_R(t)=1\ \text{for }0\leq t\leq R,\qquad
\chi_R(t)=0\ \text{for }t\geq 2R.
\]
Set
\[
r_R:=\chi_R r.
\]
Then $r_R\in C_0([0,\infty))$. Moreover,
\[
|r(t)-r_R(t)|\leq \varepsilon t
\qquad\text{for all }t\geq 0.
\]

By Lemma~\ref{lem:diff_convex}, for every $\delta>0$ there exist convex functions
$q_1,q_2\in C_0([0,\infty))$ such that
\[
\|r_R-(q_1-q_2)\|_\infty<\delta.
\]
The convex Szegő theorem applied to $q_1$ and $q_2$ gives
\[
\lim_{\alpha\to\infty}
\frac{1}{c(\alpha)}
\operatorname{Tr}\bigl((q_1-q_2)(T_\sigma^\alpha)\bigr)
=
\int_X(q_1-q_2)(\sigma(x))\,d\nu(x).
\]
Since $c(\alpha)=\dim\mathcal H_\alpha$, the uniform approximation implies
\[
\left|
\frac{1}{c(\alpha)}
\operatorname{Tr}\bigl((r_R-q_1+q_2)(T_\sigma^\alpha)\bigr)
\right|
\leq \delta,
\]
and
\[
\left|
\int_X
\bigl(r_R-(q_1-q_2)\bigr)(\sigma(x))\,d\nu(x)
\right|
\leq \delta\,\nu(X).
\]
Letting $\delta\to0$, we obtain
\[
\lim_{\alpha\to\infty}
\frac{1}{c(\alpha)}
\operatorname{Tr}\bigl(r_R(T_\sigma^\alpha)\bigr)
=
\int_X r_R(\sigma(x))\,d\nu(x).
\]

Finally, since $T_\sigma^\alpha\geq0$,
\[
\left|
\frac{1}{c(\alpha)}
\operatorname{Tr}\bigl((r-r_R)(T_\sigma^\alpha)\bigr)
\right|
\leq
\varepsilon\,
\frac{1}{c(\alpha)}
\operatorname{Tr}(T_\sigma^\alpha)
=
\varepsilon\int_X\sigma\,d\nu,
\]
while
\[
\left|
\int_X(r-r_R)(\sigma(x))\,d\nu(x)
\right|
\leq
\varepsilon\int_X\sigma\,d\nu.
\]
Letting first $\alpha\to\infty$ and then $\varepsilon\to0$ proves
\[
\lim_{\alpha\to\infty}
\frac{1}{c(\alpha)}
\operatorname{Tr}\bigl(r(T_\sigma^\alpha)\bigr)
=
\int_X r(\sigma(x))\,d\nu(x).
\]
Adding the linear term $\beta x$ completes the proof.
\end{proof}


\subsection{Determinant-type Szeg\H{o} theorem}

Although Theorem~\ref{thm:general_szego} is formulated for continuous functions, its
conclusion extends to the logarithmic function by a simple
regularization argument.  This yields the determinant formulation of
the first Szeg\H{o} theorem.

\begin{thm}\label{thm:determinant_szego}
Suppose that the assumptions of Theorem~\ref{thm:general_szego} hold and that
$c(\alpha)=\dim H_\alpha$.  Let $\sigma\in L^1(X,\nu)$ be non-negative. Then
\[
 \lim_{\alpha\to\infty}
 \frac{1}{c(\alpha)}
 \text{tr}\left( \log\left(T_\sigma^\alpha\right)\right)
 =
 \int_X \log\sigma(x)\,d\nu(x),
\]
where both sides are understood in the extended real sense. In
particular, the right-hand side belongs to $[-\infty,\infty)$.
\end{thm}

\begin{proof}
For $\varepsilon>0$, let
\[
 \psi_\varepsilon(x)=\log(x+\varepsilon).
\]
Since $\psi_\varepsilon$ is continuous on $[0,\infty)$ and
$\psi_\varepsilon(x)/x\to0$ as $x\to\infty$, Theorem~4.3 gives
\[
 \lim_{\alpha\to\infty}
 \frac{1}{c(\alpha)}
 \tr\log(T_\sigma^\alpha+\varepsilon I)
 =
 \int_X\log(\sigma(x)+\varepsilon)\,d\nu(x).
\]
Since
\[
 T_\sigma^\alpha\leq T_\sigma^\alpha+\varepsilon I,
\]
the monotonicity of the logarithm gives
\[
 \frac{1}{c(\alpha)}\text{tr}\left(\log\left( T_\sigma^\alpha\right)\right)
 \leq
 \frac{1}{c(\alpha)}
 \tr\log(T_\sigma^\alpha+\varepsilon I).
\]
Consequently,
\[
 \limsup_{\alpha\to\infty}
 \frac{1}{c(\alpha)}\text{tr}\left(\log\left( T_\sigma^\alpha\right)\right)
 \leq
 \int_X\log(\sigma+\varepsilon)\,d\nu.
\]
Letting $\varepsilon\downarrow0$ and an application of monotone convergence yields
\[
 \limsup_{\alpha\to\infty}
 \frac{1}{c(\alpha)}\text{tr}\left(\log\left( T_\sigma^\alpha\right)\right)
 \leq
 \int_X\log\sigma\,d\nu,
\]
where the integral is understood in the extended sense.

For the reverse inequality, fix $\varepsilon>0$ and define
\[
\phi_\varepsilon(t)
=
\begin{cases}
-\log(t+\varepsilon), & t\geq 0,\\[2mm]
-\log(\varepsilon)-\dfrac{t}{\varepsilon}, & t<0.
\end{cases}
\]
Then $\phi_\varepsilon$ is convex on $\mathbb{R}$ and agrees with
$-\log(t+\varepsilon)$ on $[0,\infty)$.

Applying the upper Berezin--Lieb inequality to
$\phi_\varepsilon$ gives
\[
\operatorname{Tr}\bigl(-\log(T_\sigma^\alpha+\varepsilon I)\bigr)
\leq
\int_X-\log(\sigma(x)+\varepsilon)\,d\nu_\alpha(x).
\]
Equivalently,
\[
\frac{1}{c(\alpha)}
\operatorname{Tr}\log(T_\sigma^\alpha+\varepsilon I)
\geq
\int_X\log(\sigma(x)+\varepsilon)\,d\nu(x).
\]
Letting $\varepsilon\downarrow0$ and using monotone convergence gives
\[
\frac{1}{c(\alpha)}
\operatorname{Tr}\log(T_\sigma^\alpha)
\geq
\int_X\log(\sigma(x))\,d\nu(x),
\]
where both sides are interpreted in the extended real sense.
\end{proof}

\subsection{Eigenvalue counting consequences} We are in the setting of Theorem~\ref{thm:general_szego}. So $X$ is compact, $(X,\nu)$ is a finite measure space, $\mathcal{H}_\alpha$ are all finite-dimensional  with 
$\dim \mathcal{H}_\alpha=c_\alpha$ for all $\alpha$,
and $T_\sigma^\alpha$ are Toeplitz operators with real-valued $\sigma\in L^1(X,\nu)$ such that either $\sigma$ bounded or $\sigma$ non-negative, and eigenvalues
\(\lambda_1^\alpha\ge\cdots\ge\lambda_{d_\alpha}^\alpha\), where $d_{\alpha}$ counts multiplicity. 

\begin{cor}
For any $t\in\mathbb R$ such that
\(\nu(\sigma=t)=0\),
\[
\lim_{\alpha\to\infty} \frac{1}{\dim \mathcal{H}_\alpha} 
\#\{j:\lambda_j^\alpha>t\}
=
\nu(\{x\in X:\sigma(x)>t\}).
\]
\end{cor}

\begin{proof}
Define the empirical eigenvalue measures
\[
\mu_\alpha = \frac{1}{d_\alpha} \sum_{j=1}^{d_\alpha} \delta_{\lambda_j^\alpha}.
\]
The abstract Szegő theorem gives
\(\mu_\alpha \Rightarrow \sigma_* \nu\) weakly.
By the Portmanteau theorem, for any $t$ with $\nu(\sigma=t)=0$,
\[
\lim_{\alpha\to\infty} \mu_\alpha((t,\infty)) = (\sigma_* \nu)((t,\infty)) = \nu(\sigma>t),
\]
which is exactly the statement.
\end{proof}

\subsection{Compact Abelian Groups}

The classical Szeg\H{o} limit theorem on $\mathbb{T}$ extends naturally to general compact metrizable Abelian groups. This setting was first considered by B\'edos \cite{bedos1996folner} and motivated two problems posed by N. Nikolski and A. Pushnitski \cite{nikolski2020szeg} in the context of infinite-dimensional tori. These problems were essentially resolved in a recent work by Guo, Li, and Zhou \cite{guo2024first}. Here, we present an improvement of B\'edos' results to \emph{non-negative integrable symbols} on compact Abelian groups. The locally compact, non-compact case will be treated in Section~5.

Let $G$ be a compact Abelian group with normalized Haar measure $\mu$, and let $\hat{G}$ denote its dual group. Let $\{\Gamma_N\}_{N\in\mathbb{N}} \subset \hat{G}$ be a sequence of finite subsets satisfying the F\o lner condition:
\[
\forall \eta \in \hat{G}, \qquad \lim_{N\to\infty} \frac{|(\eta \cdot \Gamma_N) \triangle \Gamma_N|}{|\Gamma_N|} = 0.
\]
The metrizability assumption is exactly what permits the sequential formulation: $\widehat G$ is then countable and admits a F{\o}lner sequence. For nonmetrizable compact Abelian groups, the corresponding formulation uses F{\o}lner nets rather than sequences.
It is well known that this condition is equivalent to requiring the limit to hold for all compact subsets $K \subset \hat{G}$.

For each $N$, define the finite-dimensional subspace 
\[
\mathcal{K}_N := \mathrm{span}\{\xi : \xi \in \Gamma_N\} \subset L^2(G),
\] 
whose elements are linear combinations of characters in $\Gamma_N$. These characters form an orthonormal basis of $\mathcal{K}_N$, which is a reproducing kernel Hilbert space with kernel
\[
K^N(x,y) = \sum_{\xi \in \Gamma_N} \xi(x-y), \qquad x,y \in G.
\]
The norm of the kernel at any point satisfies
\[
\|K^N_x\|^2 = K^N(x,x) = |\Gamma_N|.
\]

Endow $G$ with the scaled Radon measures $\mu_N := |\Gamma_N| \, \mu$. Then the normalized reproducing kernels 
\[
k^N_x := \frac{K^N_x}{\|K^N_x\|}, \quad x \in G,
\]
form a continuous Parseval frame for $\mathcal{K}_N$, since for any $f \in \mathcal{K}_N$ we have
\[
\|f\|^2 = \int_G |f(x)|^2 \, d\mu(x) = \int_G |\langle f, k^N_x \rangle|^2 \, d\mu_N(x).
\]
In particular, Assumption (I) of our general Szeg\H{o} framework is automatically satisfied, with $c(N)=|\Gamma_N|=\dim \mathcal{K}_N$.

Moreover, the kernel satisfies
\[
|\Gamma_N| \, |\langle k^N_x, k^N_y \rangle|^2 = \frac{1}{|\Gamma_N|} \left| \sum_{\xi \in \Gamma_N} \xi(x-y) \right|^2,
\]
which is a standard approximate identity on $G$. This follows, for example, from Lemma 3.4 of \cite{coifman1973operators}, and verifies Assumption (II). Hence, the general Szeg\H{o} theorem applies in this setting. Besides Bedos's theorem we also obtain Theorem~\ref{thm:compact_ab_group}, which is a  generalization of~\cite{guo2024first} for all compact Abelian groups.

\begin{proof}[Proof of Theorem~\ref{thm:compact_ab_group}]
The space $\mathcal{K}_N$ and normalized kernels $k^N_x$ satisfy Assumptions (I) and (II) of our general Szeg\H{o} theorem. Moreover, the normalized kernel $|\langle k^N_x, k^N_y \rangle|^2$ acts as an approximate identity on $G$, ensuring that the Berezin transform of $T_\sigma^N$ converges in measure to $\sigma$. The result then follows directly by applying Theorem~\ref{thm:general_szego} to $\mathcal{K}_N$ and $\sigma$.
\end{proof}

As a further consequence of the abstract determinant theorem, we obtain
the determinant formulation of the first Szeg\H{o} theorem on compact
Abelian groups.  In contrast with formulations requiring the symbol to
be bounded away from zero, the following result allows arbitrary
non-negative integrable symbols.

\begin{cor}\label{DetGroup}
Let $G$ be a compact metrizable Abelian group with normalized Haar measure $\mu$,
and let $\{\Gamma_N\}$ be a F{\o}lner sequence as above.  Let
$\sigma\in L^1(G,\mu)$ be non-negative, and let
$T_\sigma^N:K_N\to K_N$ be the corresponding Toeplitz operator.  Then
\[
    \lim_{N\to\infty}
    \frac{1}{|\Gamma_N|}
    \text{tr}\left(\log\left( T_\sigma^N\right)\right)
    =
    \int_G \log\sigma(x)\,d\mu(x),
\]
where the equality is understood in the extended real sense.  In
particular, the right-hand side is allowed to equal $-\infty$.
\end{cor}

\begin{proof}
As shown above, the spaces $K_N$, together with their normalized
reproducing kernels, satisfy Assumptions~(I) and~(II), with
\[
    c(N)=|\Gamma_N|=\dim K_N.
\]
The result therefore follows directly from Theorem~\ref{thm:determinant_szego}.
\end{proof}

This extends the determinant formulation to non-negative integrable
symbols without requiring the symbol to be bounded away from zero.
Together with Theorem~\ref{thm:compact_ab_group}, it shows that both the trace and determinant
forms of the first Szeg\H{o} limit theorem follow from the same
continuous-frame argument.

\section{Abstract Szeg\H{o} Limit Theorem for Locally Compact Spaces}

Throughout this section let $X$ be locally compact and let
$\sigma,\eta \in L^1(X,\nu)$ be real-valued.
As usual we assume that the normalized measures satisfy

\[
d\nu_\alpha = c(\alpha)\, d\nu,
\qquad\text{(Assumption (I))}.
\]

We also assume that the Berezin transforms satisfy

\[
\widetilde{\sigma}^{\alpha} \to \sigma,
\qquad
\widetilde{\eta}^{\alpha} \to \eta
\quad\text{in $\nu$-measure as } \alpha\to\infty,
\]
which we know is true under (Assumption (II)). 

Finally, we assume that the Toeplitz operators $T_\sigma^\alpha$ and $T_\eta^\alpha$ commute for all $\alpha$.

\subsection{Convex Szeg\H{o} Limit}

\begin{thm}\label{thm:convex_szego_lc} Let $\sigma,\eta\in L^1(X,\nu)\cap L^2(X,\nu)$ be real-valued with $\eta$ non-negative.
Assume that the Toeplitz operators $T_\sigma^\alpha$ and
$T_\eta^\alpha$ commute for all $\alpha$. 
\leavevmode
\begin{itemize}
\item[(i)] Suppose $\sigma \ge 0$ and $\psi \in C([0,\infty))$
is convex such that the limit
\[
\lim_{x\to\infty} \frac{\psi(x)}{x}
\]
exists. Then
\[
\lim_{\alpha\to\infty}
\frac{1}{c(\alpha)}
\tr\!\big(T_\eta^\alpha \psi(T_\sigma^\alpha)\big)
=
\int_X \eta(x)\psi(\sigma(x))\,d\nu(x).
\]

\item[(ii)] If $\sigma \in L^\infty(X,\nu)$ and
$\psi$ is convex and continuous on
$[\mathrm{ess\,inf}\,\sigma,\mathrm{ess\,sup}\,\sigma]$,
then the same limit holds.
\end{itemize}
\end{thm}

\begin{proof} We prove (i); part (ii) follows by a simpler argument.

\medskip
\noindent
\textbf{Step 1: The linear case.}

For $\psi(x)=x$ we have
\[
\frac{1}{c(\alpha)}
\tr\!\big(T_\eta^\alpha T_\sigma^\alpha\big)
=
\int_X \eta(x)\widetilde{\sigma}^\alpha(x)\, d\nu(x),
\]
using Assumption (I).

By Proposition~\ref{Lp},
\[
\widetilde{\sigma}^{\,\alpha}\longrightarrow\sigma
\qquad\text{in }L^2(X,\nu).
\]
Therefore, by the Cauchy--Schwarz inequality,
\[
\left|
\int_X\eta(x)
\bigl(\widetilde{\sigma}^{\,\alpha}(x)-\sigma(x)\bigr)
\,d\nu(x)
\right|
\leq
\|\eta\|_{L^2(X,\nu)}
\left\|
\widetilde{\sigma}^{\,\alpha}-\sigma
\right\|_{L^2(X,\nu)}
\longrightarrow 0.
\]
Consequently,
\[
\frac{1}{c(\alpha)}
\operatorname{Tr}\!\left(
T_\eta^\alpha T_\sigma^\alpha
\right)
\longrightarrow
\int_X\eta(x)\sigma(x)\,d\nu(x).
\]

\medskip
\noindent
\textbf{Step 2: Reduction to sublinear growth.}

Write
\[
\psi(x)=\beta x + \psi_0(x),
\qquad
\text{where } \frac{\psi_0(x)}{x}\to 0.
\]
By linearity it suffices to treat the case
$\psi(x)/x\to 0$.

\medskip
\noindent
\textbf{Step 3: Berezin--Lieb inequality.}

For convex $\psi$ the weighted Berezin--Lieb inequality yields
\[
\int_X \eta \psi(\widetilde{\sigma}^\alpha)\, d\nu_\alpha
\le
\tr\!\big(T_\eta^\alpha \psi(T_\sigma^\alpha)\big)
\le
\int_X \widetilde{\eta}^\alpha \psi(\sigma)\, d\nu_\alpha.
\]

Dividing by $c(\alpha)$ and using Assumption (I),
\[
\int_X \eta \psi(\widetilde{\sigma}^\alpha)\, d\nu
\le
\frac{1}{c(\alpha)}
\tr\!\big(T_\eta^\alpha \psi(T_\sigma^\alpha)\big)
\le
\int_X \widetilde{\eta}^\alpha \psi(\sigma)\, d\nu.
\]

\medskip

\noindent\textbf{Step 4: Convergence of both bounds.}

Since $\psi_0(x)/x\to0$, there exists a constant $C>0$ such that
\[
|\psi_0(t)|\leq t+C,
\qquad t\geq0.
\]

By Proposition~\ref{Lp},
\[
\widetilde{\sigma}^{\,\alpha}\to\sigma
\quad\text{and}\quad
\widetilde{\eta}^{\,\alpha}\to\eta
\qquad\text{in }L^2(X,\nu).
\]
Therefore,
\[
\left\|
\eta\bigl(\widetilde{\sigma}^{\,\alpha}-\sigma\bigr)
\right\|_{L^1}
\leq
\|\eta\|_{L^2}
\left\|
\widetilde{\sigma}^{\,\alpha}-\sigma
\right\|_{L^2}
\longrightarrow0.
\]
Thus
\[
\eta\,\widetilde{\sigma}^{\,\alpha}
\longrightarrow
\eta\sigma
\qquad\text{in }L^1(X,\nu).
\]
Moreover,
\[
\bigl|\eta\,\psi_0(\widetilde{\sigma}^{\,\alpha})\bigr|
\leq
\bigl|\eta\,\widetilde{\sigma}^{\,\alpha}\bigr|+C\bigl|\eta\bigr|.
\]
The right-hand side is a uniformly integrable family. Since
\[
\eta\,\psi_0(\widetilde{\sigma}^{\,\alpha})
\longrightarrow
\eta\,\psi_0(\sigma)
\]
in measure, Vitali's convergence theorem gives
\[
\int_X
\eta\,\psi_0(\widetilde{\sigma}^{\,\alpha})\,d\nu
\longrightarrow
\int_X\eta\,\psi_0(\sigma)\,d\nu.
\]

Similarly,
\[
\left\|
\sigma\bigl(\widetilde{\eta}^{\,\alpha}-\eta\bigr)
\right\|_{L^1}
\leq
\|\sigma\|_{L^2}
\left\|
\widetilde{\eta}^{\,\alpha}-\eta
\right\|_{L^2}
\longrightarrow0.
\]
Hence
\[
\widetilde{\eta}^{\,\alpha}\sigma
\longrightarrow
\eta\sigma
\qquad\text{in }L^1(X,\nu),
\]
and
\[
\bigl|
\widetilde{\eta}^{\,\alpha}\psi_0(\sigma)
\bigr|
\leq
\bigl|\widetilde{\eta}^{\,\alpha}\sigma\bigr|
+
C\bigl|\widetilde{\eta}^{\,\alpha}\bigr|.
\]
Again, the right-hand side is uniformly integrable. Since
\[
\widetilde{\eta}^{\,\alpha}\psi_0(\sigma)
\longrightarrow
\eta\psi_0(\sigma)
\]
in measure, Vitali's theorem yields
\[
\int_X
\widetilde{\eta}^{\,\alpha}\psi_0(\sigma)\,d\nu
\longrightarrow
\int_X\eta\psi_0(\sigma)\,d\nu.
\]

The squeeze theorem now proves the claim.
\end{proof}

\subsection{General Continuous Functions}

\begin{thm}\label{thm:general_szego_lc}
Let $\sigma,\eta\in L^1(X,\nu)\cap L^2(X,\nu)$ be real-valued with $\eta$ non-negative. Assume that the Toeplitz operators $T_\sigma^\alpha$ and $T_\eta^\alpha$ commute for all $\alpha$.

\begin{enumerate}
\item[(i)] If $\sigma$ is bounded and
$\psi$ is continuous on
$[\mathrm{ess\,inf}\,\sigma,\mathrm{ess\,sup}\,\sigma]$,
then
\[
\lim_{\alpha\to\infty}
\frac{1}{c(\alpha)}
\tr\!\big(T_\eta^\alpha \psi(T_\sigma^\alpha)\big)
=
\int_X \eta(x)\psi(\sigma(x))\, d\nu(x).
\]

\item[(ii)] If $\sigma\ge0$ and
$\psi\in C([0,\infty))$ such that the limit
$\lim_{x\to\infty}\psi(x)/x$ exists,
then the same limit holds.
\end{enumerate}
\end{thm}

\begin{proof}
We use Lemma~\ref{lem:diff_convex} together with
Theorem~\ref{thm:convex_szego_lc}.

\medskip
\noindent
\textbf{Case (i).}
Assume $\sigma$ is bounded and $\psi$ is continuous on
$[\mathrm{ess\,inf}\,\sigma,\mathrm{ess\,sup}\,\sigma]$.

By Lemma~\ref{lem:diff_convex}(i), for every $\varepsilon>0$
there exist convex functions $\psi_1,\psi_2$ on this interval such that
\[
\|\psi-(\psi_1-\psi_2)\|_\infty < \varepsilon.
\]
Set $\varphi := \psi-(\psi_1-\psi_2)$.
Then $\|\varphi\|_\infty < \varepsilon$.

By functional calculus,
\[
\psi(T_\sigma^\alpha)
=
(\psi_1-\psi_2)(T_\sigma^\alpha)
+
\varphi(T_\sigma^\alpha).
\]
Hence
\[
\frac{1}{c(\alpha)}
\tr\!\big(T_\eta^\alpha \psi(T_\sigma^\alpha)\big)
=
\frac{1}{c(\alpha)}
\tr\!\big(T_\eta^\alpha (\psi_1-\psi_2)(T_\sigma^\alpha)\big)
+
R_\alpha,
\]
where
\[
R_\alpha
=
\frac{1}{c(\alpha)}
\tr\!\big(T_\eta^\alpha \varphi(T_\sigma^\alpha)\big).
\]

To estimate $R_\alpha$ we use the trace inequality
\[
|\tr(AB)| \le \|B\|\, \tr(|A|),
\]
valid for trace-class $A$ and bounded $B$.
Since $\|\varphi(T_\sigma^\alpha)\|
\le \|\varphi\|_\infty < \varepsilon$,
we obtain
\[
|R_\alpha|
\le
\varepsilon \,
\frac{1}{c(\alpha)} \tr(|T_\eta^\alpha|).
\]

Since $\eta\in L^1(X,\nu)$, write
\[
\eta=\eta_+-\eta_-,
\qquad
\eta_+,\eta_-\geq0.
\]
Then
\[
\begin{aligned}
\|T_\eta^\alpha\|_{\mathfrak S_1}
&\leq
\|T_{\eta_+}^\alpha\|_{\mathfrak S_1}
+
\|T_{\eta_-}^\alpha\|_{\mathfrak S_1}\\
&=
\operatorname{Tr}(T_{\eta_+}^\alpha)
+
\operatorname{Tr}(T_{\eta_-}^\alpha)\\
&=
\int_X|\eta(x)|\,d\nu_\alpha(x)\\
&=
c(\alpha)\int_X|\eta(x)|\,d\nu(x).
\end{aligned}
\]
Consequently,
\[
|R_\alpha|
\leq
\varepsilon\,
\frac{\|T_\eta^\alpha\|_{\mathfrak S_1}}{c(\alpha)}
\leq
\varepsilon\int_X|\eta(x)|\,d\nu(x).
\]

On the other hand, Theorem~\ref{thm:convex_szego_lc} applied to the convex functions
$\psi_1$ and $\psi_2$ gives
\[
\frac{1}{c(\alpha)}
\operatorname{Tr}\!\left(
T_\eta^\alpha
(\psi_1-\psi_2)(T_\sigma^\alpha)
\right)
\longrightarrow
\int_X
\eta(x)(\psi_1-\psi_2)(\sigma(x))
\,d\nu(x).
\]

Since
\[
\left\|
\psi-(\psi_1-\psi_2)
\right\|_\infty<\varepsilon,
\]
we have
\[
\left|
\int_X
\eta(x)
\bigl((\psi_1-\psi_2)(\sigma(x))-\psi(\sigma(x))\bigr)
\,d\nu(x)
\right|
\leq
\varepsilon\int_X|\eta(x)|\,d\nu(x).
\]

Combining this estimate with the bound for $R_\alpha$, we obtain
\[
\begin{aligned}
&\limsup_{\alpha\to\infty}
\left|
\frac{1}{c(\alpha)}
\operatorname{Tr}\!\left(
T_\eta^\alpha\psi(T_\sigma^\alpha)
\right)
-
\int_X\eta(x)\psi(\sigma(x))\,d\nu(x)
\right|
\\
&\qquad\leq
2\varepsilon\int_X|\eta(x)|\,d\nu(x).
\end{aligned}
\]
Since $\varepsilon>0$ is arbitrary, it follows that
\[
\lim_{\alpha\to\infty}
\frac{1}{c(\alpha)}
\operatorname{Tr}\!\left(
T_\eta^\alpha\psi(T_\sigma^\alpha)
\right)
=
\int_X\eta(x)\psi(\sigma(x))\,d\nu(x).
\]

\medskip
\noindent
\textbf{Case (ii).}

Let
\[
\beta:=\lim_{x\to\infty}\frac{\psi(x)}{x},
\qquad
r(x):=\psi(x)-\beta x.
\]
Then
\[
\frac{r(x)}{x}\longrightarrow0
\qquad (x\to\infty).
\]

Fix $\varepsilon>0$. Choose $R>0$ such that
\[
|r(t)|\leq\varepsilon t
\qquad\text{for all }t\geq R.
\]
Let $\chi_R\in C_c([0,\infty))$ satisfy
\[
0\leq\chi_R\leq1,\qquad
\chi_R(t)=1\ \text{for }0\leq t\leq R,\qquad
\chi_R(t)=0\ \text{for }t\geq2R.
\]
Set
\[
r_R:=\chi_R r.
\]
Then $r_R\in C_0([0,\infty))$ and
\[
|r(t)-r_R(t)|\leq\varepsilon t
\qquad\text{for all }t\geq0.
\]

By Lemma~\ref{lem:diff_convex}(ii), for every $\delta>0$ there exist convex functions
$q_1,q_2\in C_0([0,\infty))$ such that, with
$q:=q_1-q_2$,
\[
\|r_R-q\|_\infty<\delta.
\]
Theorem~\ref{thm:convex_szego_lc} applied to $q_1$ and $q_2$ gives
\[
\frac{1}{c(\alpha)}
\operatorname{Tr}\!\left(
T_\eta^\alpha q(T_\sigma^\alpha)
\right)
\longrightarrow
\int_X\eta(x)q(\sigma(x))\,d\nu(x).
\]

Since $\eta\geq0$,
\[
\frac{1}{c(\alpha)}
\left|
\operatorname{Tr}\!\left(
T_\eta^\alpha(r_R-q)(T_\sigma^\alpha)
\right)
\right|
\leq
\delta\,
\frac{\operatorname{Tr}(T_\eta^\alpha)}{c(\alpha)}
=
\delta\int_X\eta\,d\nu,
\]
and
\[
\left|
\int_X\eta(x)(r_R-q)(\sigma(x))\,d\nu(x)
\right|
\leq
\delta\int_X\eta\,d\nu.
\]
Letting $\delta\to0$, we obtain
\[
\lim_{\alpha\to\infty}
\frac{1}{c(\alpha)}
\operatorname{Tr}\!\left(
T_\eta^\alpha r_R(T_\sigma^\alpha)
\right)
=
\int_X\eta(x)r_R(\sigma(x))\,d\nu(x).
\]

Because $T_\eta^\alpha$ and $T_\sigma^\alpha$ commute and
$T_\eta^\alpha\geq0$, they have a common orthonormal eigenbasis.
Using $|r-r_R|\leq\varepsilon\,\mathrm{id}$ on $[0,\infty)$, we obtain
\[
\left|
\operatorname{Tr}\!\left(
T_\eta^\alpha(r-r_R)(T_\sigma^\alpha)
\right)
\right|
\leq
\varepsilon\,
\operatorname{Tr}\!\left(
T_\eta^\alpha T_\sigma^\alpha
\right).
\]
Hence, by the linear case of Theorem~\ref{thm:convex_szego_lc},
\[
\limsup_{\alpha\to\infty}
\frac{1}{c(\alpha)}
\left|
\operatorname{Tr}\!\left(
T_\eta^\alpha(r-r_R)(T_\sigma^\alpha)
\right)
\right|
\leq
\varepsilon\int_X\eta(x)\sigma(x)\,d\nu(x).
\]
Similarly,
\[
\left|
\int_X\eta(x)(r-r_R)(\sigma(x))\,d\nu(x)
\right|
\leq
\varepsilon\int_X\eta(x)\sigma(x)\,d\nu(x).
\]
Letting first $\alpha\to\infty$, then $\delta\to0$, and finally
$\varepsilon\to0$, gives
\[
\lim_{\alpha\to\infty}
\frac{1}{c(\alpha)}
\operatorname{Tr}\!\left(
T_\eta^\alpha r(T_\sigma^\alpha)
\right)
=
\int_X\eta(x)r(\sigma(x))\,d\nu(x).
\]

Finally, applying the linear case to the term $\beta x$ yields
\[
\lim_{\alpha\to\infty}
\frac{1}{c(\alpha)}
\operatorname{Tr}\!\left(
T_\eta^\alpha\psi(T_\sigma^\alpha)
\right)
=
\int_X\eta(x)\psi(\sigma(x))\,d\nu(x).
\]


\end{proof}

\begin{prop}\label{diagonal} Assume that $\sigma\in L^1(X,\nu)$ is real-valued and bounded.
Let $\psi$ be continuous on a compact interval containing
$\operatorname{ess\,ran}(\sigma)$ and $0$. Then
\[
\lim_{\alpha\to\infty}
\frac{1}{c(\alpha)}
\operatorname{Tr}
\left(
T_\sigma^\alpha\psi(T_\sigma^\alpha)
\right)
=
\int_X \sigma(x)\psi(\sigma(x))\,d\nu(x).
\]
\end{prop}

\begin{proof}
Set
\[
\varphi(t)=t\psi(t).
\]
Then $\varphi(0)=0$. Fix $\varepsilon>0$. By the Weierstrass
approximation theorem, choose a polynomial $p$ such that
\[
\|\psi-p\|_\infty<\varepsilon
\]
on the chosen compact interval, and set
\[
q(t)=t p(t).
\]
Then
\[
|\varphi(t)-q(t)|\leq \varepsilon |t|.
\]

Since $q(0)=0$, write $q=q_1-q_2$ with $q_1$ and $q_2$
convex and satisfying $q_1(0)=q_2(0)=0$. For example, if
\[
M=\inf q'',
\]
then one may take
\[
q_1=q-\frac{M}{2}t^2,
\qquad
q_2=-\frac{M}{2}t^2
\]
when $M<0$, and $q_1=q$, $q_2=0$ when $M\geq0$.

Applying Proposition~\ref{BL-classical} to $q_1$ and $q_2$, and using
$\widetilde{\sigma}^{\alpha}\to\sigma$ in $L^1(X,\nu)$, gives
\[
\lim_{\alpha\to\infty}
\frac{1}{c(\alpha)}
\operatorname{Tr}q(T_\sigma^\alpha)
=
\int_X q(\sigma(x))\,d\nu(x).
\]

Moreover,
\[
\left|\varphi(t)-q(t)\right|\leq\varepsilon |t|,
\]
so functional calculus gives
\[
\left|
\frac{1}{c(\alpha)}
\operatorname{Tr}
\bigl(\varphi-q\bigr)(T_\sigma^\alpha)
\right|
\leq
\varepsilon\,
\frac{\operatorname{Tr}|T_\sigma^\alpha|}{c(\alpha)}
\leq
\varepsilon\int_X|\sigma|\,d\nu.
\]
The same estimate gives
\[
\left|
\int_X(\varphi-q)(\sigma(x))\,d\nu(x)
\right|
\leq
\varepsilon\int_X|\sigma|\,d\nu.
\]
Letting $\varepsilon\to0$ proves the result.
\end{proof}

\subsection{Eigenvalue distribution: the locally compact case} All the assumptions in this section are still in place. Let $\{\lambda_j^\alpha\}$ denote the eigenvalues of $T_\sigma^\alpha$.

\begin{cor}\label{eigenvalue} Suppose all the assumptions from Theorem~\ref{thm:general_szego_lc} hold. For every $t\in\mathbb R$ such that
\[
\nu(\sigma=t)=0,
\]
and every non-negative $\eta\in L^1(X,\nu)\cap L^2(X,\nu)$,
\[
\lim_{\alpha\to\infty}
\frac{1}{c(\alpha)}
\sum_{\lambda_j^\alpha>t}
\langle T_\eta^\alpha e_j^\alpha,e_j^\alpha\rangle
=
\int_{\{\sigma>t\}} \eta(x)\, d\nu(x),
\]
where $\{e_j^\alpha\}$ is an orthonormal eigenbasis of $T_\sigma^\alpha$.
\end{cor}

\begin{proof} Choose continuous functions $\psi_{1,n},\psi_{2,n}$ such that
\[
0\leq\psi_{1,n}
\leq\mathbf 1_{(t,\infty)}
\leq\psi_{2,n}\leq1,
\]
\[
\psi_{1,n}(x)=\psi_{2,n}(x)
=\mathbf 1_{(t,\infty)}(x)
\qquad\text{for }x\notin
[t-1/n,t+1/n],
\]
and
\[
\psi_{j,n}(x)\longrightarrow
\mathbf 1_{(t,\infty)}(x)
\qquad\text{for every }x\neq t.
\]
It is clear that such functions satisfy the assumptions of Theorem~\ref{thm:general_szego_lc} (ii). Furthermore, note that for each $n$
\begin{align*}
    \text{tr}\left(T_{\eta}^{\alpha}\psi_{1,n}\left(T_{\sigma}^{\alpha}\right)\right)\leq\text{tr}\left(T_{\eta}^{\alpha} 1_{(t,\infty)}\left(T_{\sigma}^{\alpha}\right)\right)\leq\text{tr}\left(T_{\eta}^{\alpha}\psi_{2,n}\left(T_{\sigma}^{\alpha}\right)\right).
\end{align*}
Consider the lower bound. First, note that 
\begin{align*}
    \frac{1}{c(\alpha)}\text{tr}\left(T_{\eta}^{\alpha}\psi_{1,n}\left(T_{\sigma}^{\alpha}\right)\right)-\int_{X}\eta\psi_{1,n}(\sigma)d\nu\leq\frac{1}{c(\alpha)}\text{tr}\left(T_{\eta}^{\alpha} 1_{(t,\infty)}\left(T_{\sigma}^{\alpha}\right)\right)-\int_{X}\eta\psi_{1,n}(\sigma)d\nu.
\end{align*}
Taking the limit as $\alpha\to\infty$ and applying Theorem~\ref{thm:general_szego_lc} (ii), we obtain 
\begin{align*}
    0\leq\liminf_{\alpha\to\infty}\frac{1}{c(\alpha)}\text{tr}\left(T_{\eta}^{\alpha} 1_{(t,\infty)}\left(T_{\sigma}^{\alpha}\right)\right)-\int_{X}\eta\psi_{1,n}(\sigma)d\nu.
\end{align*}
Hence, 
\begin{align*}
    \int_{X}\eta\psi_{1,n}(\sigma)d\nu\leq\liminf_{\alpha\to\infty}\frac{1}{c(\alpha)}\text{tr}\left(T_{\eta}^{\alpha} 1_{(t,\infty)}\left(T_{\sigma}^{\alpha}\right)\right).
\end{align*}
Since $\nu(\sigma=t)=0$,
\[
\psi_{j,n}(\sigma(x))
\longrightarrow
\mathbf 1_{\{\sigma>t\}}(x)
\qquad\text{for $\nu$-almost every }x.
\]
Moreover,
\[
\left|\eta(x)\psi_{j,n}(\sigma(x))\right|
\leq |\eta(x)|.
\]
Therefore, the Dominated Convergence Theorem gives
\[
\int_X\eta(x)\psi_{j,n}(\sigma(x))\,d\nu(x)
\longrightarrow
\int_{\{\sigma>t\}}\eta(x)\,d\nu(x).
\]

For the upper bound, we similarly have
\[
\frac{1}{c(\alpha)}
\operatorname{Tr}\!\left(
T_\eta^\alpha
\mathbf{1}_{(t,\infty)}(T_\sigma^\alpha)
\right)
\le
\frac{1}{c(\alpha)}
\operatorname{Tr}\!\left(
T_\eta^\alpha
\psi_{2,n}(T_\sigma^\alpha)
\right).
\]
Therefore, by Theorem~\ref{thm:general_szego_lc},
\[
\limsup_{\alpha\to\infty}
\frac{1}{c(\alpha)}
\operatorname{Tr}\!\left(
T_\eta^\alpha
\mathbf{1}_{(t,\infty)}(T_\sigma^\alpha)
\right)
\le
\int_X \eta(x)\psi_{2,n}(\sigma(x))\,d\nu(x).
\]
Letting \(n\to\infty\) and applying the dominated convergence theorem gives
\[
\limsup_{\alpha\to\infty}
\frac{1}{c(\alpha)}
\operatorname{Tr}\!\left(
T_\eta^\alpha
\mathbf{1}_{(t,\infty)}(T_\sigma^\alpha)
\right)
\le
\int_{\{\sigma>t\}}\eta(x)\,d\nu(x).
\]
Together with the lower bound, this proves the claimed limit.

\end{proof}
\begin{cor}
\label{cor:eig-count-sigma}
Let $\sigma\in L^1(X,\nu)$ be real-valued. Assume, in addition, $\sigma$ is bounded or non-negative. 
Let $\{\lambda_j^{(\alpha)}\}$ denote the eigenvalues of $T^\alpha_\sigma$ in
decreasing order.  If $t>0$ satisfies $\nu(\{\sigma=t\})=0$, then
\[
\lim_{\alpha\to\infty}
\frac{1}{c(\alpha)}\,\#\{j:\lambda_j^{(\alpha)}>t\}
=
\nu(\sigma>t).
\]
\end{cor}

\begin{proof} 
First suppose that $\sigma\geq0$. We use the following consequence of
Proposition~\ref{BLweighted}: for every $\varphi\in C_0([0,\infty))$,
\begin{equation}\label{prop2.2}
\frac{1}{c(\alpha)}
\operatorname{Tr}\!\left(
T_\sigma^\alpha\varphi(T_\sigma^\alpha)
\right)
\longrightarrow
\int_X\sigma(x)\varphi(\sigma(x))\,d\nu(x).
\end{equation}

Indeed, by Lemma~\ref{lem:diff_convex}(ii), choose convex functions
$q_1,q_2\in C_0([0,\infty))$ such that
\[
\|\varphi-(q_1-q_2)\|_\infty<\varepsilon.
\]
Applying Proposition~\ref{BLweighted} with $\eta=\sigma$ and $\psi=q_j$ gives
\[
\int_X\sigma q_j(\widetilde{\sigma}^{\,\alpha})\,d\nu
\leq
\frac{1}{c(\alpha)}
\operatorname{Tr}\!\left(
T_\sigma^\alpha q_j(T_\sigma^\alpha)
\right)
\leq
\int_X\widetilde{\sigma}^{\,\alpha}q_j(\sigma)\,d\nu.
\]
Since $\widetilde{\sigma}^{\,\alpha}\to\sigma$ in $L^1(X,\nu)$,
the two bounds converge to
\[
\int_X\sigma(x)q_j(\sigma(x))\,d\nu(x).
\]
The uniform approximation then gives \eqref{prop2.2}, because
\[
\left|
\frac{1}{c(\alpha)}
\operatorname{Tr}\!\left(
T_\sigma^\alpha
\bigl(\varphi-(q_1-q_2)\bigr)(T_\sigma^\alpha)
\right)
\right|
\leq
\varepsilon\int_X\sigma\,d\nu.
\]

Choose continuous functions $g_{1,n},g_{2,n}$ such that
\[
0\leq g_{1,n}
\leq \mathbf 1_{(t,\infty)}
\leq g_{2,n}\leq1,
\]
and
\[
g_{j,n}(x)\longrightarrow \mathbf 1_{(t,\infty)}(x)
\qquad\text{for every }x\neq t.
\]
For sufficiently large $n$, define
\[
\varphi_{j,n}(x)
=
\begin{cases}
\dfrac{g_{j,n}(x)}{x},&x>0,\\[2mm]
0,&x=0.
\end{cases}
\]
Since $t>0$, the functions $\varphi_{j,n}$ belong to
$C_0([0,\infty))$. For every eigenvalue $\lambda_j^\alpha\geq0$,
\[
\lambda_j^\alpha\varphi_{1,n}(\lambda_j^\alpha)
\leq
\mathbf 1_{(t,\infty)}(\lambda_j^\alpha)
\leq
\lambda_j^\alpha\varphi_{2,n}(\lambda_j^\alpha).
\]
Therefore,
\[
\frac{1}{c(\alpha)}
\sum_j\lambda_j^\alpha
\varphi_{1,n}(\lambda_j^\alpha)
\leq
\frac{1}{c(\alpha)}
\#\{j:\lambda_j^\alpha>t\}
\leq
\frac{1}{c(\alpha)}
\sum_j\lambda_j^\alpha
\varphi_{2,n}(\lambda_j^\alpha).
\]

For all sufficiently large $n$, we may choose $g_{j,n}$ to vanish on
$[0,t/2]$. Hence, for $j=1,2$,
\[
0\leq \sigma(x)\phi_{j,n}(\sigma(x))
 =g_{j,n}(\sigma(x))
 \leq \mathbf{1}_{\{\sigma>t/2\}}(x)
 \leq \frac{2}{t}\sigma(x).
\]

Using \eqref{prop2.2} and then letting $n\to\infty$, the Dominated Convergence Theorem gives
\[
\lim_{\alpha\to\infty}
\frac{1}{c(\alpha)}
\#\{j:\lambda_j^\alpha>t\}
=
\nu(\{\sigma>t\}),
\]
because $\nu(\{\sigma=t\})=0$.

If \(\sigma\) is bounded but not necessarily non-negative, use
Proposition~\ref{diagonal} instead of the previous corollary. For sufficiently large \(n\),
choose the functions \(g_{j,n}\) so that they vanish in a neighborhood
of \(0\), and define
\[
\phi_{j,n}(x)=
\begin{cases}
g_{j,n}(x)/x, & x>0,\\
0, & x\le 0.
\end{cases}
\]
Then \(\phi_{j,n}\) is continuous on every compact interval containing
the essential range of \(\sigma\) and \(0\). Proposition 5.3 gives
\[
\frac{1}{c(\alpha)}
\sum_j \lambda_j^{(\alpha)}
\phi_{j,n}(\lambda_j^{(\alpha)})
\longrightarrow
\int_X \sigma(x)\phi_{j,n}(\sigma(x))\,d\nu(x).
\]
Since \(t>0\), the same cutoff inequalities and dominated convergence
argument yield
\[
\frac{1}{c(\alpha)}
\#\{j:\lambda_j^{(\alpha)}>t\}
\longrightarrow
\nu(\{\sigma>t\}).
\]
\end{proof}

\subsection{Weighted Bergman Spaces on the Unit Ball $\mathbb{B}^n$}

Let $n\ge 1$ and let $\mathbb{B}^n\subset\mathbb{C}^n$ denote the unit ball.  
For $\alpha>-1$ we consider the weighted Bergman space
\[
A^2_\alpha(\mathbb{B}^n)
=
\left\{
f\in\mathcal{O}(\mathbb{B}^n):
\;
\|f\|^2_{A^2_\alpha}
:=
c_{\alpha,n}
\int_{\mathbb{B}^n} |f(z)|^2 (1-|z|^2)^\alpha\, dV(z)
<\infty
\right\},
\]
where $dV$ is Lebesgue measure on $\mathbb{C}^n$ and
$c_{\alpha,n}$ is the normalizing constant
\[
c_{\alpha,n}
=
\frac{\Gamma(n+\alpha+1)}{\pi^n\,\Gamma(\alpha+1)}
\]
so that $\|1\|_{A^2_\alpha}=1$.
It is well known that $A^2_\alpha(\mathbb{B}^n)$ is a reproducing kernel Hilbert space with reproducing kernel
\begin{equation}
K^\alpha(z,w)
=
\frac{1}{(1-\langle z,w\rangle)^{n+1+\alpha}},
\qquad z,w\in\mathbb{B}^n.
\label{eq:ball-kernel}
\end{equation}
In particular,
\[
\|K^\alpha_z\|^2_{A^2_\alpha}
=
K^\alpha(z,z)
=
(1-|z|^2)^{-(n+1+\alpha)}.
\]

\medskip

\noindent
Thus the normalized reproducing kernels are
\[
k^\alpha_z(w)
=
\frac{K^\alpha(w,z)}{\sqrt{K^\alpha(z,z)}},
\qquad 
z,w\in\mathbb{B}^n.
\]
A direct computation using \eqref{eq:ball-kernel} yields
\[
\left|
\langle k^\alpha_z, k^\alpha_w\rangle_{A^2_\alpha}
\right|
=
\left(
\frac{(1-|z|^2)(1-|w|^2)}
{|1-\langle z,w\rangle|^2}
\right)^{\frac{n+1+\alpha}{2}}.
\]
Recall that the Bergman distance $\beta(z,w)$ satisfies
\[
\cosh^2\!\left(\frac{\beta(z,w)}{2}\right)
=
\frac{|1-\langle z,w\rangle|^2}{(1-|z|^2)(1-|w|^2)},
\]
so that
\[
\frac{(1-|z|^2)(1-|w|^2)}{|1-\langle z,w\rangle|^2}
=
\operatorname{sech}^2\!\left(\frac{\beta(z,w)}{2}\right).
\]
Therefore
\begin{equation}
\left|
\langle k^\alpha_z, k^\alpha_w\rangle_{A^2_\alpha}
\right|^2
=
\operatorname{sech}^{\,2(n+1+\alpha)}\!\left(\frac{\beta(z,w)}{2}\right).
\label{eq:kernel-decay}
\end{equation}

\medskip

\noindent
The associated measures in the abstract framework are
\[
d\nu_\alpha(z)
=
\|K^\alpha_z\|^2_{A^2_\alpha}\,
c_{\alpha,n}(1-|z|^2)^\alpha\, dV(z)
=
c_{\alpha,n}\,(1-|z|^2)^{-n-1}\, dV(z).
\]
Since $c_{\alpha,n}$ depends on $\alpha$ while 
$(1-|z|^2)^{-n-1} dV(z)$ does not,  
Assumption~\textup{(I)} holds with
\[
c(\alpha)=c_{\alpha,n},
\qquad 
d\nu(z)=(1-|z|^2)^{-n-1}dV(z),
\]
independent of $\alpha$.

\medskip

\noindent
We now verify Assumption~\textup{(II)}.
Let $\beta(\cdot,\cdot)$ denote the Bergman distance on $\mathbb{B}^n$, and let
\[
\mathbb{B}^n(z,R)
=
\{w\in\mathbb{B}^n : \beta(w,z)<R\}
\]
be the Bergman ball of radius $R>0$.
Using the invariance of both the Bergman metric and the measure $d\nu_\alpha$
under the automorphism group $\operatorname{Aut}(\mathbb{B}^n)$, it suffices to treat the case $z=0$.
In this case $\mathbb{B}^n(0,R)$ is the Euclidean ball of radius $s=\tanh (R/2)$.

For $z=0$, the inner product formula \eqref{eq:kernel-decay} reduces to
\[
\left|\langle k^\alpha_0, k^\alpha_w\rangle\right|^2
=
(1-|w|^2)^{n+1+\alpha}.
\]
Hence
\begin{align*}
&\int_{\mathbb{B}^n(0,R)^c}
\left|\left\langle k_0^\alpha,k_w^\alpha\right\rangle\right|^2
\,d\nu_\alpha(w) \\
&\qquad=
c_{\alpha,n}\int_{|w|\geq s}(1-|w|^2)^\alpha\,dV(w) \\
&\qquad=
c_{\alpha,n}\sigma_{2n-1}
\int_s^1 (1-r^2)^\alpha r^{2n-1}\,dr \\
&\qquad=
\frac{c_{\alpha,n}\sigma_{2n-1}}{2}
\int_0^{1-s^2}u^\alpha(1-u)^{n-1}\,du \\
&\qquad\leq
\frac{c_{\alpha,n}\sigma_{2n-1}}{2(\alpha+1)}
(1-s^2)^{\alpha+1}.
\end{align*}
Since
\[
c_{\alpha,n}
=
\frac{\Gamma(n+\alpha+1)}{\pi^n\Gamma(\alpha+1)}
=
\frac{1}{\pi^n}\prod_{j=1}^n(\alpha+j)
=
O(\alpha^n),
\]
and $0<1-s^2<1$, the final expression tends to $0$ as
$\alpha\to\infty$. 
Thus, Assumption~\textup{(II)} follows.

\medskip

\noindent
We therefore obtain Theorem~\ref{Sze2} for weighted Bergman spaces on the unit ball $\mathbb{B}^n$
as a direct consequence of the general Szeg\H{o}-type limit theorems, Theorem~\ref{thm:general_szego_lc} and
Proposition~\ref{diagonal}.

\subsection{Locally Compact Abelian Groups}

Let $G$ be a locally compact Abelian group with Haar measure $m$ and dual group $\widehat{G}$ equipped with Haar measure $\widehat{m}$.  
Let $\{F_N\}_{N\ge1}\subset \widehat{G}$ be a F\o lner sequence, i.e.
\[
\frac{\widehat{m}((F_N+\gamma)\triangle F_N)}{\widehat{m}(F_N)}
\longrightarrow 0
\qquad \text{for every } \gamma\in\widehat{G}.
\]

Define the Paley–Wiener type spaces
\[
\mathcal{H}_N
=
\{f\in L^2(G): \operatorname{supp}(\widehat{f})\subset F_N\}.
\]
Each $\mathcal{H}_N$ is a reproducing kernel Hilbert space with kernel
\[
K_N(x,y)
=
\widehat{\chi_{F_N}}(x-y),
\]
where $\chi_{F_N}$ denotes the characteristic function of $F_N$ and $\widehat{\cdot}$ denotes the inverse Fourier transform on $\widehat{G}$.

Since
\[
K_N(x,x)
=
\widehat{\chi_{F_N}}(0)
=
\widehat{m}(F_N),
\]
the normalization constant is
\[
c(N)=\widehat{m}(F_N).
\]

Thus Assumption (I),
\[
d\nu_N = c(N)\, dm,
\]
is satisfied.

Moreover, the normalized kernel satisfies
\[
c(N)\, |\langle k_x^N,k_y^N\rangle|^2
=
\frac{1}{\widehat{m}(F_N)}
\big|\widehat{\chi_{F_N}}(x-y)\big|^2,
\]
which forms an approximate identity on $G$ by the F\o lner property.
In particular, the Berezin transforms converge in measure for
$L^1$-symbols.

Therefore, Theorem~\ref{thm:general_szego_lc} and Proposition~\ref{diagonal}
apply in this setting and yields Theorem 1.3.

\medskip







\medskip







\section{Further Applications}

In this final section we illustrate the scope of our abstract Szeg\H{o} limit theorems by presenting several further applications.  
These include classical results of Landau--Widom for Paley--Wiener spaces, 
a recent Wiener-type lemma of Jaming--Kellay--Pérez on locally compact Abelian groups, 
a Szeg\H{o}-type theorem for Gabor (time--frequency) localization operators originally due to Feichtinger and Nowak, and finally, the necessary density condition for sampling in the Bargmann-Fock space, originally due to Lindholm.  
All of these follow naturally from Theorem~\ref{thm:general_szego_lc} under the appropriate geometric choices of $X$, $\nu$, $c(\alpha)$, and the normalized kernel family $\{k^\alpha_x\}$.

\subsection{Gabor Localization Operators}

Fix any window
\[
g\in L^2(\mathbb R),
\qquad
\|g\|_{L^2(\mathbb R)}=1.
\]
We use the time--frequency shifts
\[
\pi(q,p)f(t)=e^{2\pi i pt}f(t-q).
\]
For $h>0$, let
\[
D_hf(t)=h^{-1/4}f(t/\sqrt h),
\qquad
g_h=D_hg,
\]
and define the semiclassical shifts by
\[
\pi_h(x,\xi)f(t)
=
e^{2\pi i \xi t/h}f(t-x).
\]
Then
\[
\pi_h(x,\xi)
=
D_h\pi(x/\sqrt h,\xi/\sqrt h)D_h^{-1}.
\]

The short-time Fourier transform associated with $g_h$ is
\[
V_{g_h}f(x,\xi)
=
\langle f,\pi_h(x,\xi)g_h\rangle_{L^2(\mathbb R)}.
\]
Let $\mathcal H_h$ denote the range of $V_{g_h}$, equipped with the
inner product for which $V_{g_h}$ is unitary.

The normalized reproducing kernels are
\[
k_z^h(w)
=
\left\langle
\pi_h(z)g_h,\pi_h(w)g_h
\right\rangle_{L^2(\mathbb R)},
\qquad z,w\in\mathbb R^2.
\]
Define
\[
\gamma(\zeta)
=
\left|
\left\langle g,\pi(\zeta)g\right\rangle_{L^2(\mathbb R)}
\right|^2.
\]
Then
\[
\left|
\left\langle k_z^h,k_w^h\right\rangle_{\mathcal H_h}
\right|^2
=
\gamma\left(\frac{w-z}{\sqrt h}\right).
\]
Moreover, the Moyal identity gives
\[
\int_{\mathbb R^2}\gamma(\zeta)\,d\zeta=1.
\]
Consequently, for every $R>0$,
\[
\int_{\{|w-z|>R\}}
\left|
\left\langle k_z^h,k_w^h\right\rangle_{\mathcal H_h}
\right|^2
h^{-1}\,dw
=
\int_{\{|\zeta|>R/\sqrt h\}}\gamma(\zeta)\,d\zeta
\longrightarrow 0
\]
as $h\to0$. Thus the frame measures are
\[
d\nu_h(z)=h^{-1}\,dz,
\qquad
c(h)=h^{-1}.
\]

For a non-negative symbol $\sigma$, define
\[
T_h(\sigma)f
=
\int_{\mathbb R^2}
\sigma(z)
\langle f,k_z^h\rangle_{\mathcal H_h}
k_z^h\,d\nu_h(z).
\]

\begin{thm}\label{FeiNow}
Let
\[
\sigma\in L^1(\mathbb R^2)\cap L^2(\mathbb R^2),
\qquad
\sigma\geq0.
\]
Let $\psi:[0,\infty)\to\mathbb R$ be continuous and suppose that
\[
\lim_{t\to\infty}\frac{\psi(t)}{t}
\]
exists and is finite. Then
\[
\lim_{h\to0}
h\,\operatorname{Tr}
\bigl(T_h(\sigma)\psi(T_h(\sigma))\bigr)
=
\int_{\mathbb R^2}
\sigma(z)\psi(\sigma(z))\,dz.
\]
\end{thm}

This result is related to the theorem of Feichtinger and Nowak
\cite{FN01}. They use the same arbitrary fixed normalized
window $g$, the shifts $\pi(q,p)$, and the dilated symbols
\[
b_R(q,p)=b\left(\frac{q}{R},\frac{p}{R}\right).
\]
Their result is
\[
\lim_{R\to\infty}
\frac{1}{R^2}
\operatorname{Tr}
\bigl(T_{b_R}\psi(T_{b_R})\bigr)
=
\int_{\mathbb R^2}
b(z)\psi(b(z))\,dz,
\]
for $b\in L^1(\mathbb R^2)$ with $0\leq b\leq1$ and
$\psi\in C([0,1])$.

Taking
\[
R=h^{-1/2}
\]
gives $R^{-2}=h$. Under the unitary dilation $D_h$, the operator
$T_h(\sigma)$ is unitarily equivalent to the Feichtinger--Nowak
operator with symbol
\[
\sigma_R(z)=\sigma(z/R).
\]
Hence the preceding semiclassical formula is precisely the
Feichtinger--Nowak limit written in semiclassical variables.

The present theorem extends their result by allowing symbols in
$L^1\cap L^2$ which need not be bounded, and by allowing continuous
functions $\psi$ on $[0,\infty)$ with linear growth at infinity.

\subsection{Paley--Wiener Spaces on $\mathbb{R}^{n}$ and Dilates of Convex Bodies}

We next consider the multidimensional Paley--Wiener setting treated in the classical
work of Landau--Widom~\cite{landau1980eigenvalue}.
Let $\Omega\subset\mathbb{R}^{n}$ be a bounded convex set with nonempty interior 
(and, for convenience, symmetric about the origin).
For $\alpha>0$ define the frequency set
\[
\alpha\Omega = \{ \alpha\xi : \xi\in\Omega\},
\]
and the corresponding Paley--Wiener space
\[
PW_{\alpha\Omega}
=
\left\{
f\in L^{2}(\mathbb{R}^{n}) : \operatorname{supp}(\widehat f)\subset \alpha\Omega
\right\}.
\]
This is a reproducing kernel Hilbert space with reproducing kernel
\[
K_{\alpha}(x,y)
=
(2\pi)^{-n}
\int_{\alpha\Omega} e^{i\xi\cdot(x-y)}\,d\xi,
\qquad x,y\in\mathbb{R}^{n}.
\]
Evaluating at the diagonal gives
\[
K_{\alpha}(x,x)
=
(2\pi)^{-n}\,|\alpha\Omega|
=
(2\pi)^{-n}\,\alpha^{n}|\Omega|.
\]
Consequently, the diagonal measures in our abstract framework take the form
\[
d\nu_{\alpha}(x)
=
\|K_{\alpha,x}\|^{2}\,dx
=
c(\alpha)\,dx,
\qquad
c(\alpha)
=
(2\pi)^{-n}\,\alpha^{n}|\Omega|,
\]
so that Assumption~\textup{(I)} is satisfied.

\medskip

The normalized kernels are
\[
k_{\alpha,x}(y)
=
\frac{K_{\alpha}(y,x)}{\sqrt{K_{\alpha}(x,x)}}
=
\frac{
(2\pi)^{-n}
\int_{\alpha\Omega} e^{i\xi\cdot(y-x)}\,d\xi}
{\sqrt{(2\pi)^{-n}\alpha^{n}|\Omega|}}.
\]

A well-known asymptotic estimate (see, for example, 
\cite[Section~3]{landau1980eigenvalue}) shows that the quantities
$|\langle k_{\alpha,x},k_{\alpha,y}\rangle|^{2}$ form an approximate identity:
for every fixed
$R>0$,
\[
\int_{\|x-y\|>R}
\bigl|\langle k_{\alpha,x},k_{\alpha,y}\rangle\bigr|^2
\,d\nu_\alpha(y)
\longrightarrow0
\qquad(\alpha\to\infty).
\]
Thus, Assumption~\textup{(II)} also holds.

\medskip

Applying Theorem~\ref{thm:general_szego_lc} to $PW_{\alpha\Omega}$ gives the following 
multidimensional generalization of the one--dimensional Paley--Wiener case:

\begin{thm}
Let $\Omega\subset\mathbb{R}^n$ be a bounded convex set with
nonempty interior, and let $PW_{\alpha\Omega}$ be as above. Let
\[
\sigma\in L^1(\mathbb{R}^n)\cap L^2(\mathbb{R}^n)
\]
be non-negative. Let $\psi\in C([0,\infty))$ satisfy $\psi(0)=0$,
and suppose that
\[
\phi(t):=
\begin{cases}
\dfrac{\psi(t)}{t},&t>0,\\[2mm]
\displaystyle\lim_{s\downarrow0}\frac{\psi(s)}{s},&t=0
\end{cases}
\]
is continuous on $[0,\infty)$. Assume also that
\[
\lim_{t\to\infty}\frac{\psi(t)}{t}
\]
exists. Then
\[
\lim_{\alpha\to\infty}
\frac{1}{c(\alpha)}
\operatorname{Tr}\psi\bigl(T_\alpha(\sigma)\bigr)
=
\int_{\mathbb{R}^n}\psi(\sigma(x))\,dx,
\]
where
\[
c(\alpha)=(2\pi)^{-n}\alpha^n|\Omega|.
\]
\end{thm}

\begin{proof} Set $\eta=\sigma$ in Theorem~\ref{thm:general_szego_lc} and apply it to the function
$\phi(x)=\psi(x)/x$. Since $\phi$ is bounded and
continuous, it satisfies the growth condition in Theorem~\ref{thm:general_szego_lc}.
Moreover,
\[
T_\alpha(\sigma)\phi\bigl(T_\alpha(\sigma)\bigr)
=
\psi\bigl(T_\alpha(\sigma)\bigr)
\]
by functional calculus, and
\[
\sigma(x)\phi(\sigma(x))=\psi(\sigma(x)).
\]
The claimed formula follows.
    
\end{proof}

This recovers the leading term in the Landau--Widom asymptotic formula 
for time--band limiting operators~\cite{landau1980eigenvalue}.
Their celebrated second-order correction term, 
involving the geometry of the boundary $\partial\Omega$,
lies beyond the reach of our current abstract method.
Nonetheless, the first-order term follows immediately from the
general Szeg\H{o} theorem in a conceptually simple and unified way.

\begin{cor}
\label{cor:slepian}
Let $\Omega=[-1,1]\subset\mathbb{R}$ and let $PW_{\alpha\Omega}$ be the Paley--Wiener
space of $L^{2}$-functions whose Fourier transforms are supported in $[-\alpha,\alpha]$.
Let $\sigma=\chi_{[-1,1]}$ be the indicator function of the time interval $[-1,1]$,
and let $T_{\alpha}(\sigma)$ be the associated time--frequency localization operator
\[
T_{\alpha}(\sigma) f 
=
\int_{-1}^{1} f(x)\, K_{\alpha}(x,\cdot) \, dx,
\]
where $K_{\alpha}$ is the reproducing kernel of $PW_{\alpha\Omega}$.
Let $\lambda^{(\alpha)}_{1}\ge \lambda^{(\alpha)}_{2}\ge \dots$ 
denote the eigenvalues of $T_{\alpha}(\sigma)$.

Then for every $t\in(0,1)$ such that $|\{x:\sigma(x)=t\}|=0$, 
\[
\lim_{\alpha\to\infty}
\frac{\pi}{\alpha}\,
\#\{j:\lambda^{(\alpha)}_{j}>t\}
=
|\{x\in[-1,1]:\sigma(x)>t\}|
=
2.
\]

Equivalently,
\[
\#\{j:\lambda^{(\alpha)}_{j}>t\}
=
\frac{\alpha}{\pi}\left|\{x\in[-1,1]:\sigma(x)>t\}\right|
+o(\alpha)
=
\frac{2\alpha}{\pi}+o(\alpha).
\]
\end{cor}

In particular, the number of eigenvalues of $T_{\alpha}(\sigma)$
that lie near $1$ grows like $\frac{2\alpha}{\pi}$.
These eigenfunctions are the classical prolate spheroidal wave functions (Slepian functions),
and this recovers the sharp Slepian--Pollak--Landau concentration law.

\subsection{A Wiener-Type Limit on LCA Groups}
Let $G$ be a locally compact Abelian group with dual $\widehat{G}$, and let $\{F_N\}$ be a Følner sequence in $\widehat{G}$.  
Section~5.5 shows that each Paley--Wiener type space 
\[
H_N = \{f\in L^{2}(G): \mathrm{supp}(\widehat{f})\subset F_N\}
\]
fits into our abstract setting, with $c(N)=\widehat{m}(F_N)$ and normalized kernels forming an approximate identity. For \(x\neq y\),
\[
\langle k_x^N,k_y^N\rangle
=
\frac{1}{\widehat m(F_N)}
\int_{F_N}\xi(x-y)\,d\widehat m(\xi).
\]
Since \(x-y\neq0\), the function
\(\xi\mapsto\xi(x-y)\) is a nontrivial character of \(\widehat G\).
The Følner property therefore implies
\[
\langle k_x^N,k_y^N\rangle\longrightarrow0.
\]
Thus, Assumption (III) holds. Therefore, Proposition~\ref{lkc} applies and 
we recover the recent Wiener-type lemma of Jaming--Kellay--Pérez \cite{jaming2026wiener}:
\[
\mu(\{x\})
=
\lim_{N\to\infty}
\frac{1}{\widehat{m}(F_N)}
\int_{F_N} \xi(x)\,\widehat{\mu}(\xi)\, d\widehat{m}(\xi),
\]
for any finite complex Borel measure $\mu$.

\subsection{Sampling in the Bargmann--Fock Space}
\label{sec:FockSampling}

In this subsection we show how the classical Lindholm--Landau density condition \cite{lindholm2001sampling}
for sampling in the Bargmann--Fock space follows directly from
our locally compact abstract Szeg\H{o} theorem, through its eigenvalue distribution consequence (Corollary~\ref{cor:eig-count-sigma}).  

Recall that
\[
\Falpha
=\Bigl\{ f\in\mathcal{O}(\C): 
\|f\|_{\Falpha}^{2}
=\frac{\alpha}{\pi}
\int_{\C} |f(z)|^{2}e^{-\alpha|z|^{2}}\,dm(z)<\infty\Bigr\},
\]
The reproducing kernel is
\[
K_\alpha(z,w)=e^{\alpha z\overline w},
\]
and the normalized kernels satisfy
\[
k_z^{(\alpha)}(w)
=e^{\alpha w\overline z-\alpha|z|^2/2},
\qquad
\left|\langle k_z^{(\alpha)},k_w^{(\alpha)}\rangle\right|^2
=e^{-\alpha|z-w|^2}.
\]
Moreover,
\[
d\nu_\alpha(z)=\frac{\alpha}{\pi}\,dm(z),
\qquad
c(\alpha)=\frac{\alpha}{\pi}.
\]
Consequently,
\[
\int_{|w-z|>R}
\left|\langle k_z^{(\alpha)},k_w^{(\alpha)}\rangle\right|^2
\,d\nu_\alpha(w)
=e^{-\alpha R^2}\longrightarrow0.
\]
Thus, the Fock spaces satisfy Assumptions (I) and (II). 

The dilation operator
\[
(D_\alpha f)(z)=f(\sqrt{\alpha}\,z).
\]
is a unitary map $\calF^{2}\to\Falpha$, and a straightforward computation shows the validity of the following lemma.

\begin{lem}
\label{lem:scalingFock}
Let $\Lambda\subset\C$, and define the scaled set
\[
\Lambda_\alpha := \frac{1}{\sqrt{\alpha}}\,\Lambda.
\]
Then $\Lambda$ is sampling (resp.\ interpolating) in $\calF^{2}$ with constants
$A,B$ if and only if $\Lambda_\alpha$ is sampling (resp.\ interpolating) in
$\Falpha$ with the same constants $A,B$.
\end{lem}

For $a\in\C$ and $R>0$ set
\[
T^\alpha_{a,R}:=P_\alpha\,M_{\chi_{B(a,R)}}\,P_\alpha,
\] where $P_\alpha$ is the Bergman projection onto $\Falpha$. 
As in the Fock model, the Weyl translation
\[
(\Weyl_b f)(z)=e^{\alpha(z\overline{b}-\frac12|b|^{2})}f(z-b)
\]
is unitary on $\Falpha$, satisfies $\Weyl_b P_\alpha=P_\alpha\Weyl_b$, and

\[W_a M_{\chi_{B(0,R)}}W_a^*
=M_{\chi_{B(a,R)}},
\qquad
T_{a,R}^\alpha
=W_aT_{0,R}^\alpha W_a^*.\]
Hence

\begin{lem}[Translation covariance]
\label{lem:WeylBall}
For every $a\in\C$, $R>0$, and $\alpha>0$,
\[
T^\alpha_{a,R}
=\Weyl_{a}\,T^\alpha_{0,R}\,\Weyl_{a}^{*}.
\]
Thus, $T^\alpha_{a,R}$ has the same eigenvalues for all $a\in\C$.
\end{lem}

Let $N_\alpha^{(R)}(t)$ denote the number of eigenvalues of $T^\alpha_{a,R}$
exceeding $t\in(0,1)$; by Lemma~\ref{lem:WeylBall} this is independent of $a$.

Applying Corollary~\ref{cor:eig-count-sigma} to the bounded non-negative symbol
\[
\sigma=\chi_{B(a,R)},
\]
we obtain, for every $t\in(0,1)$,
\[
\frac{1}{c(\alpha)}
N_\alpha^{(R)}(t)
\longrightarrow
|\{z\in\mathbb C:\chi_{B(a,R)}(z)>t\}|
=
|B(a,R)|
=
\pi R^2.
\]
Since $c(\alpha)=\alpha/\pi$, this is equivalent to
\begin{equation}\label{eq:FockEigenvalueCounting}
\frac{N_\alpha^{(R)}(t)}{\alpha R^2}
\longrightarrow1.
\end{equation}


Assume that $\Lambda$ is a separated sampling set for $F^2$.
Thus, there exists $A>0$ such that
\[
A\|f\|_{F^2}^2
\leq
\sum_{\lambda\in\Lambda}
|\langle f,k_\lambda^{(1)}\rangle_{F^2}|^2,
\qquad
f\in F^2.
\]
By the dilation lemma, $\Lambda_\alpha
=\alpha^{-1/2}\Lambda$ is a sampling set for $F_\alpha^2$ with the
same lower sampling constant:
\[
A\|f\|_{F_\alpha^2}^2
\leq
\sum_{\lambda\in\Lambda_\alpha}
|\langle f,k_\lambda^{(\alpha)}\rangle_{F_\alpha^2}|^2.
\]

Fix $R>0$ and choose $0<\delta<R$. Also choose
\[
0<\rho_0<
\min\left\{\frac12\operatorname{sep}(\Lambda),R\right\},
\qquad
\rho_\alpha=\frac{\rho_0}{\sqrt{\alpha}}.
\]
Since
\[
\operatorname{sep}(\Lambda_\alpha)
=
\frac{\operatorname{sep}(\Lambda)}{\sqrt{\alpha}},
\]
the balls
\[
\{B(\lambda,\rho_\alpha):
\lambda\in\Lambda_\alpha\}
\]
are pairwise disjoint.

A standard Bargmann--Fock mean-value estimate (see \cite{Z12} Lemma 2.32) gives, for
$w\notin B(a,R)$,
\[
|\langle f,k_w^{(\alpha)}\rangle|^2
\leq
\frac{\alpha}{\pi}
\left(1-e^{-\alpha\rho_\alpha^2}\right)^{-1}
\int_{B(w,\rho_\alpha)}
|\langle f,k_u^{(\alpha)}\rangle|^2\,du.
\]
Since $\alpha\rho_\alpha^2=\rho_0^2$, summing over
$\lambda\in\Lambda_\alpha\setminus B(a,R)$ gives
\[
\begin{aligned}
\sum_{\lambda\in\Lambda_\alpha\setminus B(a,R)}
|\langle f,k_\lambda^{(\alpha)}\rangle|^2
&\leq
\left(1-e^{-\rho_0^2}\right)^{-1} \\
&\quad\times
\int_{\mathbb C\setminus B(a,R-\rho_\alpha)}
|\langle f,k_u^{(\alpha)}\rangle|^2
\frac{\alpha}{\pi}\,du.
\end{aligned}
\tag{6.2}
\]

Choose $t\in(0,1)$ sufficiently close to $1$ that
\[
\left(1-e^{-\rho_0^2}\right)^{-1}(1-t)
\leq \frac{A}{2}.
\]
Let $E_\alpha(a,R-\delta;t)$ be the span of the eigenvectors of
$T_{a,R-\delta}^{(\alpha)}$ whose eigenvalues are strictly greater
than $t$. Hence
\[
\dim E_\alpha(a,R-\delta;t)
=
N_\alpha^{(R-\delta)}(t).
\]

For all sufficiently large $\alpha$, we have $\rho_\alpha<\delta$.
If $f\in E_\alpha(a,R-\delta;t)$ and
$\|f\|_{F_\alpha^2}=1$, then
\[
\int_{B(a,R-\delta)}
|\langle f,k_u^{(\alpha)}\rangle|^2
\frac{\alpha}{\pi}\,du
=
\langle T_{a,R-\delta}^{(\alpha)}f,f\rangle
>t.
\]
Using the Parseval frame identity,
\[
\int_{\mathbb C}
|\langle f,k_u^{(\alpha)}\rangle|^2
\frac{\alpha}{\pi}\,du
=
1,
\]
and the inclusion
\[
B(a,R-\delta)\subset B(a,R-\rho_\alpha),
\]
we obtain
\[
\int_{\mathbb C\setminus B(a,R-\rho_\alpha)}
|\langle f,k_u^{(\alpha)}\rangle|^2
\frac{\alpha}{\pi}\,du
\leq 1-t.
\]
Therefore, by \emph{(6.2)},
\[
\begin{aligned}
\sum_{\lambda\in\Lambda_\alpha\cap B(a,R)}
|\langle f,k_\lambda^{(\alpha)}\rangle|^2
&\geq
A-
\left(1-e^{-\rho_0^2}\right)^{-1}(1-t) \\
&\geq \frac{A}{2}.
\end{aligned}
\]

Consequently, the restriction map
\[
\mathcal R_{\alpha,a}:
E_\alpha(a,R-\delta;t)
\longrightarrow
\ell^2\bigl(\Lambda_\alpha\cap B(a,R)\bigr),
\]
defined by
\[
\mathcal R_{\alpha,a}f
=
\bigl(
\langle f,k_\lambda^{(\alpha)}\rangle
\bigr)_{\lambda\in\Lambda_\alpha\cap B(a,R)},
\]
is injective. Hence
\[
N_\alpha^{(R-\delta)}(t)
\leq
\#\bigl(\Lambda_\alpha\cap B(a,R)\bigr),
\qquad
a\in\mathbb C.
\]

By Corollary~\ref{cor:eig-count-sigma}, applied to the symbol
$\chi_{B(a,R-\delta)}$, we have
\[
\frac{N_\alpha^{(R-\delta)}(t)}{c(\alpha)}
\longrightarrow
|B(a,R-\delta)|
=
\pi(R-\delta)^2,
\]
where
\[
c(\alpha)=\frac{\alpha}{\pi}.
\]
Thus
\[
\frac{N_\alpha^{(R-\delta)}(t)}
{\alpha(R-\delta)^2}
\longrightarrow 1.
\]
Taking the infimum over $a$ gives
\[
\liminf_{\alpha\to\infty}
\frac{
\inf_{a\in\mathbb C}
\#\bigl(\Lambda_\alpha\cap B(a,R)\bigr)
}{
\alpha(R-\delta)^2
}
\geq 1.
\]

Since $\Lambda_\alpha=\alpha^{-1/2}\Lambda$,
\[
\#\bigl(\Lambda\cap B(z,R\sqrt{\alpha})\bigr)
=
\#\bigl(\Lambda_\alpha\cap B(z/\sqrt{\alpha},R)\bigr).
\]
Therefore,
\[
\liminf_{\alpha\to\infty}
\frac{
\inf_{z\in\mathbb C}
\#\bigl(\Lambda\cap B(z,R\sqrt{\alpha})\bigr)
}{
\pi(R\sqrt{\alpha})^2
}
\geq
\frac{(R-\delta)^2}{\pi R^2}.
\]
Since $\alpha$ ranges over all positive real numbers, $R\sqrt{\alpha}$
ranges over all sufficiently large radii. Hence
\[
D^-(\Lambda)
\geq
\frac{(R-\delta)^2}{\pi R^2}.
\]
Finally, letting $\delta\downarrow0$ yields
\[
D^-(\Lambda)\geq\frac1\pi.
\]

\bibliographystyle{plain}
\bibliography{refs}
\end{document}